\documentclass[11pt]{article}
\usepackage{amssymb}
\usepackage{times}

\title{How can we recover Baire class one functions?\indent}
\author{Dominique LECOMTE}
\date{\it ~Mathematika~\rm 50 (2003), 171-198}

\newcommand{\Ana}{{\it\Sigma}^{1}_{1}}

\newcommand{\Borel}{{\it\Delta}^{1}_{1}}

\newcommand{\ana}{{\bf\Sigma}^{1}_{1}}
\newcommand{\ca}{{\bf\Pi}^{1}_{1}}
\newcommand{\borel}{{\bf\Delta}^{1}_{1}}
\newcommand{\boraone}{{\bf\Sigma}^{0}_{1}}
\newcommand{\boratwo}{{\bf\Sigma}^{0}_{2}}
\newcommand{\borathree}{{\bf\Sigma}^{0}_{3}}

\newcommand{\boraxi}{{\bf\Sigma}^{0}_{\xi}}

\newcommand{\bortwo}{{\bf\Delta}^{0}_{2}}
\newcommand{\borone}{{\bf\Delta}^{0}_{1}}
\newcommand{\borthree}{{\bf\Delta}^{0}_{3}}

\newcommand{\bormone}{{\bf\Pi}^{0}_{1}}
\newcommand{\bormtwo}{{\bf\Pi}^{0}_{2}}
\newcommand{\bormthree}{{\bf\Pi}^{0}_{3}}

\newtheorem{thm} {Theorem}
\newtheorem{defi} [thm] {Definition}
\newtheorem{cor} [thm] {Corollary}
\newtheorem{lem} [thm] {Lemma}

\newtheorem{prop} [thm] {Proposition}

\begin{document}

\maketitle

\noindent {\footnotesize {\bf Abstract.} Let $X$ and $Y$ be separable metrizable spaces, and 
$f:X\rightarrow Y$ be a function. We want to recover $f$ from its values on a small set via a simple algorithm. We show that this is possible if f is Baire class one, and in fact we get a characterization. This leads us to the study of sets of Baire class one functions and to a characterization of the
separability of the dual space of an arbitrary Banach space.}\bigskip\smallskip

\noindent\bf {\Large 1 Introduction.}\rm\bigskip

 This paper is the continuation of a study by U. B. Darji and M. J. Evans in [DE]. 
We specify the term ``simple algorithm" used in the abstract. We work in separable metrizable 
 spaces $X$ and $Y$, and $f$ is a function from $X$ into $Y$. Recall that $f$ is Baire 
 class one if the inverse image of each open set is $F_\sigma$. Assume that we only know the values 
 of $f$ on a countable dense set $D\subseteq X$. We want to recover, in a simple way, all the 
 values of $f$. For each point $x$ of $X$, we extract a subsequence of $D$ which tends to $x$. 
 Let $(s_n[x,D])_n$ be this sequence. We will say that $f$ is
$recoverable~with~respect~to~D$ if, for each $x$ in $X$, the sequence $\big(
f(s_n[x,D])\big)_n$ tends to $f(x)$. The function $f$ is $recoverable$ if there exists 
$D$ such that $f$ is recoverable with respect to $D$. Therefore, continuous functions are 
recoverable with respect to any countable dense sequence in $X$. We will show 
that results concerning recoverability depend on the way of extracting the subsequence. 
We let $D:=(x_p)$.

\begin{defi} Let $X$ be a topological space. We say that 
a basis $(W_m)$ for the topology of $X$ is a $good~basis$ if
for each open subset $U$ of $X$ and each point $x$ of $U$, there exists an integer 
$m_0$ such that, for each $m\geq m_0$, $W_m\subseteq U$ if $x\in W_m$.\end{defi}

 We show that every separable metrizable space has a good basis, 
using the embedding into the compact space $[0,1]^\omega$. In the sequel, $(W_m)$
will be a good basis of $X$, except where indicated.

\begin{defi} Let $x\!\in\! X$. The   
$path~to~x~based~on~D$ is the sequence 
$\big(s_n[x,D]\big)_{n\in\omega}$, denoted ${\cal R}\big(x,D\big)$, defined 
by induction as follows:
$$\left\{\!\!
\begin{array}{ll}
s_0[x,D] 
& \!\!\!\! := x_0\mbox{,}\cr\cr 
s_{n+1}[x,D] 
& \!\!\!\! := \left\{\!\!\!\!\!\!\!\!
\begin{array}{ll} 
& s_n[x,D]~~\mbox{if}~~x= s_n[x,D]\mbox{,}\cr\cr
& x_{\mbox{min}\Bigl\{p~/~\exists m\in\omega~\{x,x_p\}\subseteq W_m\subseteq
X\setminus\bigl\{s_0[x,D],...,s_n[x,D]\bigr\}\Bigr\}}~\mbox{otherwise.}
\end{array}
\right.
\end{array}
\right.$$\end{defi} 
 
 Now the definition of a $recoverable$ function is complete.
 
\vfill\eject
 
  In Section 2, we show the\bigskip
 
\noindent\bf Theorem 4\it\ A function $f$ is recoverable if and only if $f$ is Baire class one.\rm\bigskip
 
 In Section 3, we study the limits of U. B. Darji and M. J. Evans's result, using their 
way of extracting the subsequence. We give some possible extensions, 
and we show that we cannot extend it to any Polish space.\bigskip 

 In Section 4, we study the question of the uniformity of sequence $(x_p)$ for a set of
Baire class one functions. We consider $A\subseteq {\cal B}_1(X,Y)$, equipped with the 
pointwise convergence topology. We study the existence of a dense sequence $(x_p)$ of $X$
such that each function of $A$ is recoverable with respect to $(x_p)$ (if this happens, we 
say that $A$ is $uniformly\ recoverable$).\bigskip

 In the first part, we give some necessary conditions for uniform recoverability. We deduce
among other things from this an example of a metrizable compact space 
$A\subseteq {\cal B}_1(2^\omega,2)$ which is not uniformly recoverable.\bigskip

 In the second part, we study the link between the uniform recoverability of $A$ and the fact
that J. Bourgain's ordinal rank is bounded on $A$. J. Bourgain wondered whether his rank was
bounded on a separable compact space $A$ when $X$ is a metrizable compact space. We show among other things that, if $X$ and $A$ are Polish spaces, then this rank is bounded (this
is a partial answer to J. Bourgain's question).\bigskip 

 In the third part, we give some sufficient conditions for uniform recoverability. We study
among other things the link between uniform recoverability and $F_\sigma$ subsets with
open vertical sections of a product of two spaces.\bigskip

 In the fourth part, we give a characterization of the separability of the 
dual space of an arbitrary Banach space:\bigskip

\noindent\bf Theorem 30\it\ Let $E$ be a Banach space, $X := [B_{E^*},w^*]$, 
${A :=\{G\lceil X/G\in B_{E^{**}}\}}$, and $Y:=\mathbb{R}$. The following statements are equivalent:\smallskip

 \noindent (a) $E^*$ is separable.\smallskip
 
\noindent (b) $A$ is metrizable.\smallskip

\noindent (c) Every singleton of $A$ is $G_\delta$.\smallskip

\noindent (d) $A$ is uniformly recoverable.\rm\bigskip

 In the fifth part, we introduce a notion similar to that of
equicontinuity, the notion of an $equi\mbox{-}$ $Baire\ class\ one$ set of functions. We give several
characterizations of it, and we use it to study similar versions of Ascoli's theorems for Baire
class one functions. Finally, the study of the link between the notion of an equi-Baire class
one set of functions and uniform recoverability is made.\bigskip\smallskip

\noindent\bf {\Large 2 A characterization of Baire class one functions.}\rm\bigskip

  As mentionned in Section 1, we show the

\begin{prop} Every separable metrizable space has a good basis.\end{prop}

\vfill\eject

\noindent\bf Proof.\rm\ Let $X$ be a separable metrizable space. Then $X$ embeds 
into the compact metric space $[0,1]^\omega$, by $\phi$. So let, for $r$ integer,
$n_r$ be an integer and $(U^r_j)_{j\leq n_r}$ be a covering of $[0,1]^\omega$ made of 
open subsets of $[0,1]^\omega$ whose diameter is at most $2^{-r}$. To get $(W_m)$, it is enough to 
enumerate the sequence $\big(\phi^{-1}(U^r_j)\big)_{r\in\omega ,~j\leq n_r}$.$\hfill\square$

\begin{thm} A function $f$ is recoverable if and only if $f$ is Baire class one.\end{thm}

 In order to prove this, we first give a lemma. It is essentially identical to U. B. Darji and M. J. Evans's proof of the ``only if" direction. But we will use it later. So we give the details. Notice that it does not really depend of the way of extracting the subsequence.
 
\begin{lem} Assume that, for 
$q\in\omega$, ${\{x\in X~/~\exists n~s_n[x,D]=x_q\}}$ is an open 
subset of $X$. If $f$ is recoverable with respect to $D$, then $f$ is 
Baire class one.\end{lem}

\noindent\bf Proof.\rm\ Let $F$ be a closed subset of $Y$. We let, for $k$ 
integer, ${O_k := \{y\in Y~/~d(y,F)<2^{-k}\}}$. This defines an open subset of $Y$
containing $F$. Let us fix an integer $k$. Let $(x_{p_j})_j$ be the subsequence of 
$D$ made of the elements of $f^{-1}(O_k)$ (we may assume that it is infinite and 
enumerated in a 1-1 way). We let, for $j$ integer, ${U_j :=\{x\in X~/~\exists n~s_n[x,D]=x_{p_j}\}}$. 
This set is an open subset of $X$ by hypothesis. Let ${H_k := \bigcap_{i\in\omega} {\rm \big[}
(\bigcup_{j\geq i} U_j)\cup\{x_{p_0},...,x_{p_{i-1}}\}{\rm \big]}}$. This set is a 
$G_\delta$ subset of $X$.\bigskip

 Let $x\in f^{-1}(O_k)$ and $i$ be an integer. Then $s_n[x,D]\in f^{-1}(O_k)$ if $n$ is 
bigger than $n_0$ and there exists $j(n)$ such that $s_n[x,D] = x_{p_{j(n)}}$; thus $x\in U_{j(n)}$. 
Either there exists $n\geq n_0$ such that $j(n)\geq i$ and $x\in \bigcup_{j\geq i} U_j$, or 
$x_{p_{j(n)}}$ is $x_{p_q}$ if $n$ is big enough, with $q<i$, and $x = x_{p_q}$. 
In both cases, $x\in H_k$.\bigskip

 If $x\in H_k$, either there exists an integer $q$ such that $x = x_{p_q}$ and 
$f(x)\in O_k$, or for each integer $i$, there exists $j\geq i$ such that $x\in U_j$, 
and $\exists n~s_n[x,D]=x_{p_j}$, and thus $f(x)\in \overline{O_k}$.\bigskip

 Therefore ${f^{-1}(F)\subseteq \bigcap_{k\in\omega} f^{-1}(O_k)\subseteq
\bigcap_{k\in\omega} H_k\subseteq \bigcap_{k\in\omega} f^{-1}(\overline{O_k})\subseteq
f^{-1}(F)}$. We deduce that
$$f^{-1}(F) = \bigcap_{k\in\omega} H_k$$ 
is a $G_\delta$ subset of $X$.$\hfill\square$\bigskip

\noindent\bf Proof of Theorem 4.\rm\ In order to show the ``only if" direction, let us show that 
Lemma 5 applies. Set 
$${\cal O}(x,D,n) := \left\{\!\!\!\!\!\!
\begin{array}{ll} 
& \emptyset~~\mbox{if}~~x= s_n[x,D]\mbox{,}\cr\cr 
& W_{\mbox{min}\Big\{m~/~\{x,s_{n+1}[x,D]\}\subseteq W_m\subseteq
X\setminus\big\{s_0[x,D],...,s_n[x,D]\big\}\Big\}}~\mbox{otherwise.}
\end{array}
\right.$$
Note that ${\cal O}(x,D,n)\not=\emptyset$ if and only if 
$x\not=s_n[x,D]$. In this case ${\cal O}(x,D,n)$ is an open neighborhood of 
$x$. If $n<n'$ and ${{\cal O}(x,D,n),~{\cal O}(x,D,n')\not=\emptyset}$, 
$s_{n+1}[x,D]\in {\cal O}(x,D,n)\setminus {\cal O}(x,D,n')$, so 
${\cal O}(x,D,n)$ is distinct from ${\cal O}(x,D,n')$. As $(W_m)$ is a good 
basis, for each open neighborhood $V$ of $x$ there exists an integer $n_0$ such that 
${\cal O}(x,D,n)\! \subseteq\!  V$ if $n\! \geq\!  n_0$, and therefore 
$s_{n+1}[x,D]\! \in\!  V$. So path to $x$ based on $D$ tends to $x$.

\vfill\eject
 
 To show that $\{x\in X~/~x_q\in {\cal R}\big(x,D\big)\}$ is an open subset of 
$X$, we may assume that $q>0$ and that $x_r\not= x_q$ if $r<q$. So let $t_0\in X$ and 
$n$ be a minimal integer such that ${s_{n+1}[t_0,D] = x_q}$. Let $m$ be a minimal 
integer such that ${\{t_0,x_q\}\subseteq W_m\subseteq 
X\setminus\{s_{0}[t_0,D], ...,s_{n}[t_0,D]\}}$.
By definition of the path, $q$ is minimal such that $x_q\in W_m$. Let us show that 
if $x\in W_m$, then $x_q\in {\cal R}\big(x,D\big)$; this will be enough 
since $t_0\in W_m$. We notice that if we let 
${p_n(x) := \mbox{min}\{p\in\omega~/~x_p = s_n[x,D]\}}$, then the sequence 
$\big(p_n(x)\big)_n$ 
increases, strictly until it may be eventually constant. We have $x\in W_m$, which is a 
subset of $X\setminus \{x_0,...,x_{q-1}\}$. Thus, as the path to $x$ based on $D$ tends to 
$x$, there exists a minimal integer $n'$ such that $p_{n'+1}(x)\geq q$. Then we have  
$x_q = s_{n'+1}[x,D]\in {\cal R}\big(x,D\big)$.\bigskip

 Let us show the ``if" direction. The proof looks 
like C. Freiling and R. W. Vallin's ones in [FV]. The main 
difference is the choice of the dense sequence, which has to be 
valid in any separable metrizable space.\bigskip 

 We say that $D$ $approximates$ $F\subseteq X$ if for all 
$x\in F\setminus D$, ${\cal R}\big(x,D\big)\setminus F$ is finite. 
Let us show that if $(F_{i})$ is a sequence of closed subsets of $X$, then there is $D\subseteq X$ which approximates each $F_{i}$.\bigskip

 Consider a countable dense sequence of $X$, and also a countable dense sequence 
of each finite intersection of the $F_{i}$'s. Put this together, to get a countable dense sequence $(q_{i})$ of $X$. This 
countable dense set is the set $D$ we are looking for. But we've got to 
describe how to order the elements of this sequence.\bigskip 

 We will construct $D$ in stages, called $D_{i}$, for each integer $i$. If $F$ is 
a finite intersection of the $F_{i}$'s and $G$ is a finite subset of $D$, we set 
$$A^F(G) := \bigcup_{m\in\omega , x\in G\setminus F, x\in W_{m}\not\subseteq 
X\setminus F}~\{q_{\mbox{min}\{i/q_{i}\in W_{m}\cap F\}}\}.$$
 Put on $2^{i}=\{\sigma_{1},\ldots ,\sigma_{2^{i}}\}$ the lexicographic ordering, and 
 let $F^\sigma := \bigcap_{j\in\sigma}~F_{j}$ for each finite subset $\sigma$ of 
 $\omega$. We set 
$$G_{0}:=\{q_{i}\}\mbox{,}~G_{k+1}:=G_{k}\cup A^{F^{\sigma_{k+1}}}(G_{k})~~~~~(\mbox{for}~k<2^{i})\mbox{,}$$
$$D_{i}:=(\bigcup_{k\leq 
2^{i}}~G_{k})\setminus (\bigcup_{l<i}~D_{l}).$$
 We order the elements of $D_{i}$ as follows. Let $\sigma^{i} (x):= \{k < 
 i/x\in F_{k}\}$. Put the elements of $D_{i}$ whose $\sigma^{i}$ is 
 $\sigma_{2^{i}}$ first (in any order). Then put the elements of $D_{i}$ whose 
 $\sigma^{i}$ is 
 $\sigma_{2^{i}-1}$. And so on, until elements of $D_{i}$ whose 
 $\sigma^{i}$ is $\sigma_{1}$.\bigskip
 
 Now let us suppose that $F_{i}$ is not approximated by $D$, with $x$ 
 as a witness and $i$ minimal. Let $y\in {\cal R}\big(x,D\big)\setminus 
 F_{i}$ such that $y$ is put into $D$ at some stage $j>i$ 
 and satisfying $x\in F_{k}\Leftrightarrow y\in F_{k}$ for each $k<i$. 
 Let $m\in\omega$ such that $x,y\in W_{m}$. 
 We have $y\notin F^{\sigma^{j} (x)}$, and $W_{m}\not\subseteq 
 X\setminus F^{\sigma^{j}  (x)}$ because $x\in F^{\sigma^{j}  (x)}$. So we can define 
 $z:=q_{\mbox{min}\{i/q_{i}\in W_{m}\cap F^{\sigma^{j}  (x)}\}}$. Then 
 $\sigma^{j} (z)>\sigma^{j} (y)$ in the lexicographic order. We have $z\in 
 A^{F^{\sigma^{j} (x)}}(\{y\})$. We conclude that $z$ is put before $y$ 
 and that $y\notin {\cal R}\big(x,D\big)$. This is the contradiction we were looking 
 for.\bigskip
 
  Now let $(Y_{p})$ be a basis for the topology of $Y$. Consider the 
  inverse images of the $Y_{p}$'s by $f$. Express each of these sets 
  as a countable union of closed sets. This gives $D$ which 
  approximates each of these closed sets. It is now clear that the set $D$ 
  is what we were looking for.$\hfill\square$
  
\vfill\eject

\noindent\bf {\Large 3 About the limits of U. B. Darji and M. J. Evans's method.}\rm\bigskip 

Let us recall the original way of extracting the subsequence. Fix a 
compatible distance $d$ on $X$.

\begin{defi} Let $x\in X$. The   
$route~to~x~based~on~D$ is the sequence 
$\big(s'_n[x,D]\big)_{n\in\omega}$, denoted ${\cal R}'\big(x,D\big)$, defined 
by induction as follows:
$$\left\{\!\!
\begin{array}{ll}
s'_0[x,D] 
& \!\!\!\! := x_0\mbox{,}\cr\cr 
s'_{n+1}[x,D] 
& \!\!\!\! := \left\{\!\!\!\!\!\!\!\!
\begin{array}{ll} 
& s'_n[x,D]~~\mbox{if}~~x= s'_n[x,D]\mbox{,}\cr\cr
& x_{\mbox{min}\bigl\{p~/~d(x,x_{p})<d(x,s'_n[x,D])\bigr\}}~\mbox{otherwise.}
\end{array}
\right.
\end{array}
\right.$$\end{defi}
  
 If $f$ is recoverable in the sense of Definition 6, we say that $f$ is $first~return~recoverable$. U. B. Darji and M. J. Evans showed the following:\bigskip
 
\noindent\bf Theorem\it\ If $f$ is first return recoverable, then $f$ is Baire 
class one. Conversely, if $f$ is Baire class one and $X$ is a compact space, then $f$ is 
first return recoverable.\rm

\begin{defi} We will say that an ultrametric space 
$(X,d)$ is $discrete$ if the following condition is satisfied: $\forall~
(d_n)_{n\in\omega}\subseteq d[X\times X]~~~~\big[(\forall~n\in\omega~~d_{n+1}<d_n)~
\Rightarrow ~(\lim_{n\rightarrow \infty} d_n = 0)\big]$.\end{defi}

 We can show the following extensions:
 
\begin{thm} Assume that $f$ is Baire one. Then $f$ is first return recoverable in 
the following cases:\smallskip

\noindent (a) $X$ is a metric space countable union of totally bounded subspaces.\smallskip
 
\noindent (b) $X$ is a discrete ultrametric space.\end{thm}

\begin{cor} Let $X$ be a metrizable separable space. 
Then there exists a compatible distance $d$ on $X$ such that for each 
$f:X\rightarrow Y$, $f$ is Baire class one if and only if $f$ is first return 
recoverable relatively to $d$.\end{cor}

 This corollary comes from the fact that we can 
find a compatible distance on $X$ making $X$ totally bounded. Now we 
will show that the notion of a first return recoverable function is a 
metric notion and not a topological one. More precisely, we will show 
that the hypothesis ``$X$ is discrete" in Theorem 8 is useful. In fact, 
we will give an example of an ultrametric 
space homeomorphic to $\omega^\omega$ in which there exists a closed subset 
whose characteristic function is not first return recoverable (notice 
that $\omega^\omega$, equipped with its usual metric, is a discrete 
ultrametric space). So the equivalence between ``$f$ is Baire class 
one " and ``$f$ is first return recoverable" depends on the choice of 
the distance. And the equivalence in Theorem 4 does not depend on the 
choice of the good basis, and is true without any restriction on $X$. 
The algorithm given in Definition 2 is given in topological terms 
only, as the notion of a Baire class one function. Furthermore,  
Definition 2 uses only countably many open subsets of $X$, namely the 
$W_{m}$'s.

\vfill\eject 

\begin{lem} Let $X$ be an ultrametric space, $t\in X$, 
$x, y\in X\setminus\{ t\}$. Then the open balls $B\big(x,d(x,t)\big[$ and 
$B\big(y,d(y,t)\big[$ are equal or disjoint.\end{lem}

\noindent\bf Proof.\rm\ Let us show that $d(x,t) = d(y,t)$ or 
$B\big(x,d(x,t)\big[\cap B\big(y,d(y,t)\big[ = \emptyset$. Let 
$$z\in B\big(x,d(x,t)\big[ \cap B\big(y,d(y,t)\big[.$$ 
If for example $d(x,t)<d(y,t)$, let $r$ be in $\big]d(x,t),d(y,t)\big[$. As $z\in B\big(x,r\big[$, 
$$B\big(x,r\big[ = B\big(z,r\big[\subseteq B\big(z,d(y,t)\big[.$$ 
As $z\in B\big(y,d(y,t)\big[$, we can write 
${B\big(z,d(y,t)\big[ = B\big(y,d(y,t)\big[\subseteq X\setminus\{ t\}}$. But this contradicts the fact that 
$t\in B\big(x,r\big[$.\bigskip

 If $B\big(x,d(x,t)\big[ \cap B\big(y,d(y,t)\big[ \not= \emptyset$, let $z$ 
be in the intersection. Then we have the sequence of equalities 
${B\big(x,d(x,t)\big[ = B\big(z,d(x,t)\big[ = B\big(z,d(y,t)\big[ = 
B\big(y,d(y,t)\big[}$.$\hfill\square$\bigskip 

 Now we introduce the counterexample. We set 
$$Z := \{Q=(q_n)_{n\in\omega}\in \mathbb{Q} ~_+^{\omega}~/~\forall~n\in\omega~~q_n<q_{n+1}~~
\mbox{and}~\mbox{lim}_{n\rightarrow\infty} q_n = +\infty\}.$$
This space is equipped with 
$$d:\left\{\!\!
\begin{array}{ll}
Z\times Z 
& \!\!\!\!\rightarrow \mathbb{R}_+ \cr\cr 
(Q,Q') 
& \!\!\!\!\mapsto \left\{\!\!\!\!\!\!\!\!
\begin{array}{ll}
& 2^{-\mbox{min}(q_{\mbox{min}\{n\in\omega/q_n\not=q'_n\}},
q'_{\mbox{min}\{n\in\omega/q_n\not=q'_n\}})}~~\mbox{if}~~Q\not= Q'\mbox{,}\cr\cr 
& 0~~\mbox{otherwise.}
\end{array}
\right.
\end{array}\right.$$

\begin{prop} The space $(Z,d)$ is an ultrametric space homeomorphic to $\omega^\omega$ and is not discrete.\end{prop} 

\noindent\bf Proof.\rm\ We set 
$W := \{f\in 2^{\mathbb{R}_+}~/~\exists~Q\in Z~~f= 1\!\! {\rm I}
_{\cup_{p\in\omega} [q_{2p},q_{2p+1}]}\}$; this space is equipped with 
the ultrametric on $2^{\mathbb{R}_+}$ defined by $\tilde d (f,g) := 
2^{-\mbox{inf}\{x\in\mathbb{R}_{+}/f(x)\not= g(x)\}}$ if $f\not= g$. 
Then the function from $Z$ into $W$ which associates $1\!\! {\rm I}
_{\cup_{p\in\omega} [q_{2p},q_{2p+1}]}$ to $Q$ is a bijective isometry. Thus, it is 
enough to show the desired properties for $W$.\bigskip

 We set 
$$D := \{f\in 2^{\mathbb{R}_+}~/~\exists~Q\in Z~~\exists~k\in\omega~~
f= 1\!\! {\rm I}_{\cup_{p<k} [q_{2p},q_{2p+1}]}~~\mbox{or}~~
f= 1\!\! {\rm I}_{\cup_{p<k} [q_{2p},q_{2p+1}]\cup [q_{2k},+\infty[}\}\mbox{,}$$
$$V := W\cup D.$$
Then $W$ and $V$ are ultrametric, viewed as subspaces of $2^{\mathbb{R}_+}$. Set  
$D$ is countable and dense in $V$, so $V$ and $W$ are separable. 

\vfill\eject

 Let $(f_p)_{p\in\omega}$ be a Cauchy sequence in $V$, and $m$ in $\omega$. There exists a 
minimal integer $N(m)$ such that, for $p,q\geq N(m)$, we have $d(f_p,f_q)\leq 2^{-m}$; that is to say $f_p(t)=f_q(t)$ for each $t<m$. We let, if $E(t)$ is the biggest 
integer less than or equal to $t$, 
$$f:\left\{\!\!
\begin{array}{ll}
\mathbb{R}_+ 
& \!\!\!\!\rightarrow 2 \cr t 
& \!\!\!\!\mapsto f_{N(E(t)+1)}(t)
\end{array}
\right.$$
If $p\geq N(m)$ and $t<m$, $N\big(E(t)+1\big)\leq N(m)$ and we have 
$${f(t) = f_{N(E(t)+1)}(t) = f_{N(m)}(t) = f_p(t)}.$$ 
Thus the sequence $(f_p)_{p\in\omega}$ tends to $f$ in $2^{\mathbb{R}_+}$. We will check that 
$f\in V$; 
this will show that $V$ is complete, thus Polish. As $W$ is a $G_\delta$ subset of 
$V$, $W$ will also be Polish.\bigskip

\noindent\bf Case~1.\rm\ $\exists~r\in\mathbb{R}_+~~\forall~t\geq r~~f(t)=0$.\bigskip

 If $p\geq N\big(E(r)+1\big)$ and $t<E(r)+1$, $f_p(t) = f(t)$; thus, the restriction of $f$ to 
$\big[0,E(r)+1\big[$ is the restriction of $1\!\! {\rm I}_{\cup_{p<k} 
[q_{2p},q_{2p+1}]}$ to this interval, and we may assume that $q_{2k-1}<E(r)+1$. 
Therefore, we have $f = 1\!\! {\rm I}_{\cup_{p<k} [q_{2p},q_{2p+1}]}$ and 
$f\in D\subseteq V$.\bigskip

\noindent\bf Case~2.\rm\ $\exists~r\in\mathbb{R}_+~~\forall~t\geq r~~f(t)=1$.\bigskip

 If $p\geq N\big(E(r)+1\big)$ and $t<E(r)+1$, then $f_p(t) = f(t)$; thus, the restriction of $f$ to 
$\big[0,E(r)+1\big[$ is the restriction of $1\!\! {\rm I}_{\cup_{p\leq k} 
[q_{2p},q_{2p+1}]}$ to this interval, and we may assume that $q_{2k}<E(r)+1$. 
Therefore, we have $f = 1\!\! {\rm I}_{\cup_{p<k} [q_{2p},q_{2p+1}]\cup 
[q_{2k},+\infty[}$ and $f\in D\subseteq V$.\bigskip

\noindent\bf Case~3.\rm\ $\forall~r\in\mathbb{R}_+~~\exists~t,u\geq r~~f(t)=0$ and $f(u)=1$.\bigskip

 Let $(r_n)_{n\in\omega}\subseteq \mathbb{R}_+$ be a strictly increasing sequence such that 
$\lim_{n\rightarrow\infty} r_n = +\infty$ and $f(r_n)=0$ for each integer $n$. If $t<E(r_n)+1$, then we have $f(t) = f_{N(E(r_n)+1)}(t)$. Thus, the restriction of $f$ to 
$[0,r_n]$ is the restriction of $1\!\! {\rm I}_{\cup_{p<k_n} [q_{2p},q_{2p+1}]}$ to 
this interval, and we may assume that $q_{2k_n-1}<r_n$. The sequence 
$(k_n)_{n\in\omega}$ is increasing, and $\lim_{n\rightarrow\infty} k_n = +\infty$ 
because $f$ is not ultimately constant. For the same reason, 
$\mbox{lim}_{n\rightarrow\infty} q_n = +\infty$. Thus 
$f = 1\!\! {\rm I}_{\cup_{p\in\omega} [q_{2p},q_{2p+1}]}\in W\subseteq 
V$.\bigskip 

 Let $f\in V$ and $m$ in $\omega$. There exists $\varepsilon\in\mathbb{Q} ~\cap]0,1[$ and 
$q\in\mathbb{Q} ~_+\cap ]m+1,+\infty[$ such that, for each $t\in ]q-\varepsilon,q+\varepsilon[$, we have 
$f(t) = 0$, or, for each $t>q-\varepsilon$, we have $f(t) = 1$. In the first case we set 
$$g:\left\{\!\!
\begin{array}{ll} 
\mathbb{R}_+ 
& \!\!\!\!\rightarrow 2 \cr t 
& \!\!\!\!\mapsto\left\{\!\!\!\!\!\!\!\!
\begin{array}{ll} 
& f(t)~~\mbox{if}~~t\notin [q-\varepsilon /2, q+\varepsilon /2]\mbox{,}\cr\cr 
& 1~\mbox{otherwise.}
\end{array}
\right.
\end{array}
\right.$$
In the second case, we set  
$$g:\left\{\!\!
\begin{array}{ll} 
\mathbb{R}_+ 
& \!\!\!\!\rightarrow 2 \cr t 
& \!\!\!\!\mapsto\left\{\!\!\!\!\!\!\!\!
\begin{array}{ll} 
& f(t)~~\mbox{if}~~t\leq q\mbox{,}\cr\cr 
& 0~\mbox{otherwise.}
\end{array}
\right.
\end{array}
\right.$$

\vfill\eject

 In both cases we have $f\not= g$, $d(f,g)\leq 2^{-m}$ and $g\in V$; this shows 
that $V$ is perfect. Moreover, as $D$ is countable and dense in $V$, $W$ is 
locally not compact. Finally, $W$ is a $0$-dimensional Polish 
space, and each of its compact subsets have empty interior; thus it is homeomorphic 
to $\omega^\omega$ (see Theorem 7.7 page 37 of 
chapter 1 in [Ke]).\bigskip

 To finish the proof, we set $f_n := 1\!\! {\rm I}_{[0,1-2^{-n-1}]\cup 
\bigcup_{p>0} [2p,2p+1]}$. We have $f_n\in W$ and $d(f_n,f_{n+1})$ is 
$2^{-1+2^{-n-1}}$, which strictly decreases to $1/2$. Thus, space $W$ is not a 
discrete ultrametric space.$\hfill\square$

\begin{thm} There exists a $\bormone (Z)$ whose characteristic function is not first return recoverable.\end{thm}

\noindent\bf Proof.\rm\ Let 
$F:= \{Q\in Z~/~\forall~n\in\omega~~n<q_n<n+1\}$, $D := (x_p)$ be a dense sequence of $Z$. Then $F$ 
is closed since fixing a finite number of coordinates is a clopen condition. We will 
show that there exists $x\in Z$ such that the sequence 
$\big(1\!\! {\rm I}_F(s_n[x,D])\big)_{n\in\omega}$ does not tend to 
$1\!\! {\rm I}_F(x)$. Let us assume that this is not the case.\bigskip

\noindent $\bullet$ We set $n_\emptyset := 0$, $B_{\emptyset} := Z$. We have $Z\setminus \{x_0\} = 
\bigcup_{j\in\omega}^{\mbox{disj.}}~ 
B\big(y_j,d(y_j,x_0)\big[$. Let 
$${n_j := \mbox{min}\{n\in\omega~/~x_n\in 
B\big(y_j,d(y_j,x_0)\big[\}}.$$
For each $x$ in $B\big(y_j,d(y_j,x_0)\big[$ we have ${B\big(y_j,d(y_j,x_0)\big[=
B\big(x_{n_j},d(x_{n_j},
x_0)\big[=B\big(x,d(x,x_0)\big[}$ and $s_1[x,D]=x_{n_j}$. Then we do this construction again. For 
$s\in\omega^{<\omega}\setminus\{\emptyset\}$, we set 
$${B_s := B\big(x_{n_s},d(x_{n_s},x_{n_{s\lceil\vert s\vert -1}})\big[}.$$
We have $B_s\setminus\{ x_{n_s}\} = \bigcup_{j\in\omega}^{\mbox{disj.}}~ B\big(y_{s^\frown j},
d(y_{s^\frown j},x_{n_s})\big[$. Let 
$${n_{s^\frown j} := \mbox{min}\{n\in\omega~/~x_n\in B\big(y_{s^\frown j},d(y_{s^\frown j},x_{n_s})\big[\}}.$$ 
For each $x$ in $B\big(y_{s^\frown j},d(y_{s^\frown j},x_{n_s})\big[$, we have 
$${B\big(y_{s^\frown j},d(y_{s^\frown j},x_{n_s})\big[\! = 
B\big(x_{n_{s^\frown j}},d(x_{n_{s^\frown j}},x_{n_s})\big[\! = 
B\big(x,d(x,x_{n_s})\big[\subseteq B_s}\mbox{,}$$ 
and also ${s_{\vert s\vert +1}[x,D] = x_{n_{s^\frown j}}}$.\bigskip

\noindent $\bullet$ For each $x$ in $Z\setminus\{x_n~/~n\in\omega\}$, there is $\alpha$ in 
$\omega^\omega$ with $x\in \bigcap_{m\in\omega} B_{\alpha\lceil m}$ and 
$s_{m}[x,D] 
= x_{n_{\alpha\lceil m}}$ for each $m$ in $\omega$. Moreover, if $x\in F$, then there exists 
 $m_0$ in $\omega$ such that $x_{n_{\alpha\lceil m}}\in F$ for each $m\geq 
 m_0$.\bigskip
 
 \noindent\bf Case~1.\rm~$\forall~s\in\omega^{<\omega}~~B_s\cap 
F=\emptyset~\mbox{or}~\exists~t\succ_{\not=}s~~B_t\cap F\not=\emptyset~\mbox{and}~
x_{n_t}\notin F$.\bigskip

 As $B_\emptyset = Z$ meets $F$ which is not empty, there exists  
$\alpha, \beta\in\omega^\omega$ such that ${0<\beta (n)<\beta (n+1)}$, 
$B_{\alpha\lceil\beta (n)}\cap F\not=\emptyset$ and $x_{n_{\alpha\lceil\beta (n)}}
\notin F$ for each $n$ in $\omega$. It is enough to show the existence of 
$x$ in $\bigcap_{m\in\omega} B_{\alpha\lceil m}$. Indeed, if we have this, we will have 
$s_m[x,D]=x_{n_{\alpha\lceil m}}$ for each $m\in\omega$. But the diameter of 
$B_{\alpha\lceil m}$ will be at most $2d\big(s_m[x,D],s_{m-1}[x,D]\big)$, thus will tend to $0$. 
As $B_{\alpha\lceil \beta (n)}$ meets $F$, we will 
deduce that $x\in F$. Thus, the sequence $\big(1\!\! {\rm I}_F(s_n[x,D])\big)_
{n\in\omega}$ will not tend to $1\!\! {\rm I}_F(x)$ since $s_{\beta 
(n)}[x,D]\notin F$.\bigskip

 As $x_{n_{\alpha\lceil m+1}}\in B_{\alpha\lceil m+1}\subseteq B_{\alpha\lceil m}$, the sequence 
$\big(d(x_{n_{\alpha\lceil m+1}},x_{n_{\alpha\lceil m}})\big)_{m\in\omega}$ is strictly 
decreasing; let $l$ be its inferior bound.

\vfill\eject\bigskip

\noindent\bf Case 1.1.\rm\ $l=0$.\bigskip

 In this case, sequence $(x_{n_{\alpha\lceil m}})_{m\in\omega}$ is a Cauchy 
sequence. Let $\Phi$ be the bijective isometry that we used at the beginning of the 
proof of Proposition 12. We set ${f_m := \Phi (x_{n_{\alpha\lceil m}})}$. Then the sequence $(f_m)_{m\in\omega}$ is a Cauchy sequence in $W\subseteq V$, thus tends 
to $f\in V$ which is complete.\bigskip

\noindent\bf Case 1.1.1.\rm\ $\exists~r\in\mathbb{R}_+~~\forall~t\geq r~~f(t)=0$.\bigskip

 We have $f = 1\!\! {\rm I}_{\cup_{p<k} [q_{2p},q_{2p+1}]}$ and, if 
$m$ is big enough, then the restriction of $f_m$ to $[0,E(r)+1[$ is the restriction of $f$ to this 
same interval, and we have $q_{2k-1}<E(r)+1$. Thus, $x_{n_{\alpha\lceil m}}$ starts with 
$<q_0,q_1,..., q_{2k-1},q^m_{2k}>$ and, if $m$ is greater 
than $p_0\geq m_0$, then $q^m_{2k}\geq 2k+1$. Let $n_0$ in $\omega$ be such that $\beta 
(n_0)>p_0$. Then $B_{\alpha\lceil\beta (n_0)}$ is disjoint from $F$ because, if $y$ 
is in $F$, then $y\notin B_{\alpha\lceil\beta (n_0)}$ since 
$d(y,x_{n_{\alpha\lceil \beta (n_0)}})\geq 2^{-y_{2k}}>2^{-2k-1}\geq 
d(x_{n_{\alpha\lceil \beta (n_0)}},x_{n_{\alpha\lceil \beta 
(n_0)-1}})$. Thus, this case is not possible.\bigskip

\noindent\bf Case 1.1.2.\rm\ $\exists~r\in\mathbb{R}_+~~\forall~t\geq r~~f(t)=1$.\bigskip

 This case is similar to case 1.1.1.\bigskip
 
\noindent\bf Case 1.1.3.\rm\ $\forall~r\in\mathbb{R}_+~~\exists~t, u\geq r~~f(t)=0~\mbox{and}~
f(u)=1$.\bigskip

 In this case, $f\in W$, thus there exists $x\in Z$ such that the sequence 
$(x_{n_{\alpha\lceil m}})_{m\in\omega}$ tends to $x$. We have  
$x\in\bigcap_{m\in\omega} B_{\alpha\lceil m}$, since otherwise we can find an integer $m'_0$ 
such that $x\notin B_{\alpha\lceil m}$ for each $m\geq m'_0$; but, as 
$x_{n_{\alpha\lceil m}}\in B_{\alpha\lceil m}$, $x$ is in $B_{\alpha\lceil m'_0}$ 
which is closed.\bigskip 

\noindent\bf Case 1.2.\rm\ $l>0$.\bigskip

 Let $r'\in\mathbb{R}$ be such that $l = 2^{-r'}$.\bigskip 
 
\noindent\bf Case 1.2.1.\rm\ $E(r')<r'$.\bigskip

 We will show that there exists $x\in \bigcap_{m\in\omega} B_{\alpha\lceil m}$. This 
will be enough. If $m$ is big enough, $d(x_{n_{\alpha\lceil m}}, 
x_{n_{\alpha\lceil m-1}})<2^{-E(r')}$. As $B_{\alpha\lceil m}$ meets $F$, let 
$y$ be in the intersection; $y$ is of the form 
$$(n+1-\varepsilon_n)_{n\in\omega}\mbox{,}$$ 
where  $\varepsilon_n\in ]0,1[$. If $m$ is big enough, then $x_{n_{\alpha\lceil m}}$ starts with 
$<1-\varepsilon_0, ...,E(r')-\varepsilon_{E(r')-1}>$. Then the term number $E(r')+1$ of  
sequence $x_{n_{\alpha\lceil m}}$ is called $x_{n_{\alpha\lceil m}}^
{E(r')}$.\bigskip 

\noindent\bf Case 1.2.1.1.\rm\ $\exists~m\in\omega~~x_{n_{\alpha\lceil m}}^{E(r')} = 
x_{n_{\alpha\lceil m+1}}^{E(r')}$.\bigskip

 In this case, as $B_{\alpha\lceil m+1}$ meets $F$, 
$x_{n_{\alpha\lceil p}}^{E(r')}$ is of the form $E(r')+1-\varepsilon_{E(r')}$ for each  
$p\geq m$. This shows that if $p$ is big enough, then $x_{n_{\alpha\lceil p}}^{E(r')+1}
\not= x_{n_{\alpha\lceil p+1}}^{E(r')+1}$. Thus we are reduced to the following case.

\vfill\eject

\noindent\bf Case 1.2.1.2.\rm\ $\forall~m\in\omega~~x_{n_{\alpha\lceil m}}^{E(r')} \not= 
x_{n_{\alpha\lceil m+1}}^{E(r')}$.\bigskip

 The sequence $(x_{n_{\alpha\lceil m}}^{E(r')})_{m\in\omega}$ is  
strictly increasing. Indeed, assume that ${x_{n_{\alpha\lceil m}}^{E(r')} >  
x_{n_{\alpha\lceil m+1}}^{E(r')}}$. Then we have  
${d(x_{n_{\alpha\lceil m}},x_{n_{\alpha\lceil m+1}}) \leq d(
x_{n_{\alpha\lceil m+1}},x_{n_{\alpha\lceil m+2}})}$, since 
${x_{n_{\alpha\lceil m+1}}^{E(r')}\not= x_{n_{\alpha\lceil m+2}}^{E(r')}}$; but this 
is absurd. Thus the sequence $(x_{n_{\alpha\lceil m}}^{E(r')})_{m\in\omega}$ is 
strictly increasing, and  
${\mbox{lim}_{m\rightarrow\infty}~x_{n_{\alpha\lceil m}}^{E(r')}= r'}$. But if the point $x$ starts 
with sequence ${<1-\varepsilon_0, ...,E(r')-\varepsilon_{E(r')-1}, q>}$, where $q\in\mathbb{Q} ~\cap ]r',+\infty [$, then $x\in\bigcap_{m\in\omega} B_{\alpha\lceil m}$ since 
$$d(x,x_{n_{\alpha\lceil m}}) = 2^{-x_{n_{\alpha\lceil m}}^{E(r')}} 
< d(x_{n_{\alpha\lceil m}},x_{n_{\alpha\lceil m-1}}) = 
2^{-x_{n_{\alpha\lceil m-1}}^{E(r')}}.$$ 
\bf Case 1.2.2.\rm\ $E(r')=r'$.\bigskip

 This case is similar to case 1.2.1; $r'-1$ plays the role that $E(r')$ played in 
the preceding case.\bigskip

\noindent\bf Case~2.\rm\ $~\exists~s\in\omega^{<\omega}~~B_s\cap 
F\not=\emptyset~\mbox{and}~\forall~t\succ_{\not=}s~~B_t\cap F\not=\emptyset~\Rightarrow~
x_{n_t}\in F$.\bigskip

 Note that, for each $x$ in $Z$ and each $q$ in $\mathbb{Q} ~_+$, there exists 
$Q$ in $Z$ such that $d(Q,x)=2^{-q}$. Indeed, there exists a minimal integer $n$ 
such that $q<x_n$, and we take $Q$ beginning with $<x_0,...,x_{n-1},q>$ if 
$x_{n-1}\not= q$; otherwise, we take $Q$ beginning with $<x_0,...,x_{n-2},x_n>$.\bigskip

 We may assume, by shifting $s$ if necessary, that $x_{n_s}\in F$ and  
$s\not=\emptyset$. Thus we have 
$$x_{n_s} = ~<1-\varepsilon^0_0, 2-\varepsilon^0_1,...>\!\mbox{,}$$ 
where $0<\varepsilon^0_i<1$. Let $j_0$ be a minimal integer such that 
${2^{\varepsilon^0_{j_0}-j_0-1}< d(x_{n_s},x_{n_{s\lceil\vert s\vert -1}})}$, and 
$$s_0 := ~<1-\varepsilon^0_0,...,j_0-\varepsilon^0_{j_0-1}>\! .$$ 
If $t\succ s$, then $x_{n_t}$ begins with $s_0$.\bigskip  

 There are $p_0$ in $\omega$ and $Q_{n_{s^\frown p_0}}$ in $F$ such that  
$d(Q_{n_{s^\frown p_0}},x_{n_s}) = 
2^{\varepsilon^0_{j_0}-j_0-1}$. Then $Q_{n_{s^\frown p_0}}$ is of the form 
$s_0^\frown <1+j_0-\varepsilon^1,2+j_0-\varepsilon^0_{j_0+1},...>\!\mbox{,}$ 
where $0<\varepsilon^1<\varepsilon^0_{j_0}$. There exists an unique integer $n_0$ such that 
$$Q_{n_{s^\frown p_0}}\in B_{s^\frown n_0} = B(Q_{n_{s^\frown p_0}}, 2^{\varepsilon^0_{j_0}-j_0-1}[.$$
As $B_{s^\frown n_0}$ meets  $F$, $x_{n_{s^\frown n_0}}\in F$. Thus the point $x_{n_{s^\frown n_0}}$ is of the form 
$${s_0^\frown <1+j_0-\varepsilon^1_{j_0},2+j_0-\varepsilon^1_{j_0+1},...>}\!\mbox{,}$$
where $0<\varepsilon^1_{j_0}<\varepsilon^0_{j_0}$. More generally, there exists 
$p_k$ in $\omega$ and $Q_{n_{s^\frown n_0^\frown ...^\frown n_{k-1}\frown p_k}}$ in 
$F$ such that $d(Q_{n_{s^\frown n_0^\frown ...^\frown n_{k-1}\frown p_k}},
x_{n_s^\frown n_0^\frown ...^\frown n_{k-1}}) = 2^{\varepsilon^k_{j_0}-j_0-1}$. Then 
$Q_{n_{s^\frown n_0^\frown ...^\frown n_{k-1}\frown p_k}}$ is of the form 
$$s_0^\frown <1+j_0-\varepsilon^{k+1},2+j_0-\varepsilon^k_{j_0+1},...>\!\mbox{,}$$ 
where $0<\varepsilon^{k+1}<\varepsilon^k_{j_0}$.

\vfill\eject

 There exists an unique integer $n_k$ 
such that $Q_{n_{s^\frown n_0^\frown ...^\frown n_{k-1}\frown p_k}}$ is in  
$B_{s^\frown n_0^\frown ...^\frown n_k}$, which is 
$B(Q_{n_{s^\frown n_0^\frown ...^\frown n_{k-1}\frown p_k}}, 2^{\varepsilon^k_{j_0}-j_0-1}[$. 
As $B_{s^\frown n_0^\frown ...^\frown n_k}$ meets 
$F$, $x_{n_{s^\frown n_0^\frown ...^\frown n_k}}$ is in $F$. Thus the point 
$x_{n_{s^\frown n_0^\frown ...^\frown n_k}}$ is of the form 
$s_0^\frown <1+j_0-\varepsilon^{k+1}_{j_0},2+j_0-\varepsilon^{k+1}_{j_0+1},...>$, where  
$0<\varepsilon^{k+1}_{j_0}<\varepsilon^k_{j_0}$.\bigskip

 We set $\gamma := <n_0,n_1,...>$ and 
$x:=s_0^\frown (j_0+1+k+\eta_k)_{k\in\omega}$, 
where $\eta_k\in \mathbb{Q} ~_+$ are chosen so that $\eta_0 :=0$ and 
$x\notin \{x_n~/~n\in\omega\}$. Then $d(x,x_{n_{s^\frown \gamma\lceil m}}) = 
2^{\varepsilon^m_{j_0}-j_0-1}$ decreases to $r>0$, and the sequence 
$(x_{n_{s^\frown \gamma\lceil m}})_{m\in\omega}$ does not tend to $x$. But  
$x\in\bigcap_{m\in\omega} B_{s^\frown \gamma\lceil m}$; thus 
$s_{m+\vert s\vert }[x,D] = 
x_{n_{s^\frown \gamma\lceil m}}$ and the sequence 
$\big(s_m[x,D]\big)_{m\in\omega}$ does not tend to $x$. But this is 
absurd.$\hfill\square$\bigskip\smallskip 

\noindent\bf {\Large 4 Study of the uniformity of the dense sequence.}\bigskip\rm

\noindent\bf (A) Necessary~conditions~for~uniform~recoverability.\rm\bigskip

 It is natural to wonder whether there exists a dense sequence $(x_p)$ of $X$ such that every
Baire class one function from $X$ into $Y$ is first return recoverable with respect to $(x_p)$.
The answer is no when $X$ is uncountable. Indeed, if we choose $x\in
X\setminus\{x_p~/~p\in\omega\}$, then 
$1\!\!{\rm I}_{\{ x\}}$ is not first return recoverable with respect to $(x_p)$. We can wonder 
whether $(x_p)$ exists for a set of Baire class one functions.\bigskip

\bf\noindent Notation\rm ~${\cal B}_1(X,Y)$ is the set of Baire class one functions 
from $X$ into $Y$, and is equipped with the pointwise convergence topology.\bigskip

 If $A$ is a subset of ${\cal B}_1(X,Y)$, then the map 
$$\phi :\left\{\!\!
\begin{array}{ll}
X\times A
& \!\!\!\!\rightarrow Y\cr (x,f) 
& \!\!\!\!\mapsto f(x)
\end{array}
\right.$$ 
has its partial functions 
$\phi(x,.)$ (respectively 
$\phi(.,f)$) continuous (respectively Baire class one). Therefore $\phi$ is Baire class two if 
$A$ is a metrizable separable space (see p 378 in [Ku]).

\begin{defi} We will say that $A\subseteq {\cal B}_1(X,Y)$
is $uniformly~recoverable$ if there exists a dense sequence $(x_p)$ of $X$ such
that every function of $A$ is recoverable with respect to $(x_p)$.\end{defi}

\begin{prop} If $A$ is uniformly recoverable and compact, then
$A$ is metrizable.\end{prop}

\noindent\bf Proof.\rm\ Let $D:=(x_p)$ be a dense sequence of $X$ such that every
function of $A$ is recoverable with respect to $D$. Let $I:A\rightarrow 
Y^\omega$ defined by $I(f) := \big(f(x_p)\big)_p$. This map is continuous by definition of 
the pointwise convergence topology. It is one-to-one 
because, if $f\not= g$ are in $A$, then there is $p\in\omega$ such that 
$f(x_p)\not=g(x_p)$. Indeed, if this were not the case, then we would have, for each $x$ in
$X$,
$$f(x) = \mbox{lim}_{n\rightarrow\infty}~ f\big( s_n[x,D]\big) = \mbox{lim}_{n\rightarrow\infty}~ g\big( s_n[x,D]\big) = g(x)$$ 
(because $f$ and $g$ are recoverable with respect to $(x_p)$).  As $A$ is compact, $I$ is a homeomorphism from $A$ onto a subset of
$Y^\omega$. Therefore, $A$ is metrizable.$\hfill\square$

\vfill\eject

\noindent\bf Example.\rm ~There are some separable compact spaces which are not metrizable, 
and whose points are $G_\delta$. For example, ``split interval" ${A := \{
f:[0,1]\rightarrow 2~/~f~\mbox{is~increasing}\}}$, viewed as a subset of ${{\cal
B}_1([0,1],2)}$, is one of them (see [T]).
$A$ is compact because it is a closed subset of $2^{[0,1]}$: 
$$f\in A\Leftrightarrow\forall~x\leq y~~f(x)=0~\mbox{or}~f(y)=1.$$ 
$A$ is separable because ${\{1\!\!\mbox{I}_{[q,1]}~/~q\in[0,1]\cap\mathbb{Q} ~\}\cup\{1\!\!{\rm I}_{]q,1]}~/~q\in[0,1]\cap\mathbb{Q} ~\}}$ is a 
countable dense subset of $A$. The family of continuous functions $\phi_x\! :\! f\!\mapsto\! f(x)$
separates  points, and for every sequence ${(x_n)\!\subseteq\! [0,1]}$, $(\phi_{x_n})_n$ does not
separates points.  Thus $A$ is not analytic and not metrizable (see Corollary 1 page 77 in
Chapter 9 of [Bo2]).  Finally, every point of $A$ is $G_\delta$; for example,
${\{1\!\!{\rm I}_{[x,1]}\} = 
\bigcap_{n\in\omega,x\geq 2^{-n}} \{f\in A~/~f(x-2^{-n})\not=1\}\cap
\{f\in A~/~f(x)\not=0\}}$. By Proposition 14, ``split interval" is not uniformly recoverable.

\begin{prop} If $A$ is uniformly recoverable 
and $Y$ is a 0-dimensional space, then $\phi$ is Baire class one.\end{prop}
 
\noindent\bf Proof.\rm\ Let $F$ be a closed subset of $Y$. We have 
$\phi(x,f)\in F\Leftrightarrow x\in f^{-1}(F)$. Remember the proof of 
Lemma 5. We 
replace the $O_k$'s by a sequence of clopen subsets of $Y$ whose intersection is $F$ (it
exists because $Y$ is a 0-dimensional space). The sequence $(x_{p_j})_j$ is finite or infinite and
enumerates in a one-to-one way the elements of ${(x_p)\cap f^{-1}(O_k)}$. We 
have ${U_j := \{t\in
X~/~x_{p_j}\in {\cal R}\big(t,D\big)\}}$ if $x_{p_j}$ exists ($U_j :=
\emptyset$ otherwise), and ${H_k:=\bigcap_{i\in\omega} [(\bigcup_{j\geq i}
U_j)\cup\{x_{p_0},...,x_{p_{i-1}}\}]}$ (in fact, between braces we have the $x_{p_j}$ that
exist, for $j<i$). So that $f^{-1}(F)=\bigcap_{k\in\omega} H_k$. The sequence $(x_{p_j})_j$ can be
defined as follows, by induction on integer $j$:
$$\begin{array}{ll}
q=p_0 
& \!\!\!\!\Leftrightarrow f(x_q)\in O_k~\mbox{and}~\forall~l<q~~f(x_l)\notin O_k\cr\cr  
q\! =\! p_{j+1} 
& \!\!\!\!\Leftrightarrow\forall~l<q~x_l\! \not=\!  x_q~\mbox{and}~\exists~r\! <\! q~~(r\! =\! p_j~\mbox{and}~f(x_q)\! \in\!  O_k~\mbox{and}~\forall~l\! \in
]r,q[\cap\omega~f(x_l)\! \notin\!  O_k)
\end{array}$$ 
We notice that the relation ``$q=p_j$" is clopen in $f$. Then we notice that 
$$x\in (\bigcup_{j\geq i} U_j)\cup\{x_{p_0},...,x_{p_{i-1}}\}$$ 
if and only if $[\exists j\! \geq\! i~~\exists q\! \in\! \omega~~q\! =\! p_j~\mbox{and}~x_q\in {\cal R}\big(x,D\big)]$ or $\exists r\! \leq\!  i$ such that $[(\forall m\! <\! r~~\exists q\! \in\! \omega~~q\! =\!  p_m)$ and 
$(\forall m\! \in\! [r,i[\cap\omega ~~\forall q\!\in\!\omega~~q\! \not=\! p_m)$ and 
$(\exists m\! <\! r~~\forall q\!\in\!\omega~~q\! \not=\! p_m~\mbox{or}~x\! =\! x_q)]$. We can deduce from this that the relation ``$x\in (\bigcup_{j\geq i}
U_j)\cup\{x_{p_0},...,x_{p_{i-1}}\}$" is $G_\delta$ in $(x,f)$; thus the relation ``$x\in H_k$" is too.
$\hfill\square$

\begin{cor} (a) There exists a continuous injection~
$I:2^\omega\rightarrow {\cal B}_1(2^\omega ,2)$ such that $I[2^\omega]$ is not uniformly
recoverable $\big($and in fact such that $\phi\notin {\cal
B}_1(2^\omega\times I[2^\omega],2)\big)$.\smallskip

\noindent (b) There exists $A\subseteq {\cal B}_1(2^\omega,2)$,
$A\approx\omega^\omega$, which is not uniformly recoverable and such that $\phi$ is in 
${\cal B}_1(2^\omega\times A ,2)$.\end{cor}

\noindent\bf Proof.\rm\ (a) Let ${{\cal S} := \{s\in
2^{<\omega}/s=\emptyset~\mbox{or}~[s\not=\emptyset~\mbox{and}~s(\vert s\vert -1) =1]\}}$ and 
$$I(\alpha ) :=\left\{\!\!
\begin{array}{ll}
2^\omega 
& \!\!\!\!\rightarrow 2\cr 
\beta 
& \!\!\!\!\mapsto \left\{\!\!\!\!\!\!\!\!
\begin{array}{ll} 
& 1~\mbox{if}~~\exists s\in {\cal S}~~[s\prec\alpha~\mbox{and}~\beta = s^\frown 0^\omega ]\mbox{,}\cr 
& 0~\mbox{otherwise.}
\end{array}
\right.
\end{array}
\right.$$
If $\alpha\lceil n=\alpha'\lceil n$ and $\alpha (n) = 1-\alpha'(n) = 0$, then ${I(\alpha
)(\alpha\lceil n^\frown 10^\omega ) = 0 = 1\! -\! I(\alpha'
)(\alpha\lceil n^\frown 10^\omega )}$. Thus $I$ is one-to-one.

\vfill\eject

 It is continuous because 
$$I(\alpha )(\beta ) = 1~\Leftrightarrow ~\left\{\!\!\!\!\!\!\!\!
\begin{array}{ll}
& \alpha\in\emptyset~~\mbox{if}~~\beta\in P_\infty:=\{\alpha\in2^\omega~/~\forall n~\exists m\geq n~~
\alpha (m)=1\}\mbox{,}\cr\cr 
& \alpha\in N_s~~\mbox{if}~~\beta= s^\frown 0^\omega~~\mbox{and}~~s\in {\cal S}. 
\end{array}
\right.$$
Moreover, $\{\beta\in2^\omega~/~I(\alpha )(\beta ) = 1\} = \{0^\omegaÊ\}\cup \bigcup_{n/\alpha
(n)=1} ~\{\alpha\lceil (n+1)^\frown 0^\omega\}\in D_2(\boraone)(2^\omega )$, thus ${I[2^\omega
]\subseteq {\cal B}_1(2^\omega ,2)}$. Let us argue by contradiction. We have 
$$\phi^{-1}(\{ 0\}) \equiv (P_\infty\times2^\omega )\cup
(\bigcup_{s\in {\cal S}}~\{s^\frown 0^\omega\}\times\check
N_s)=\bigcup_n F_n\in\boratwo(2^\omega\times2^\omega).$$
The diagonal of $P_\infty$ is a subset of $\phi^{-1}(\{ 0\})$, so there exists an integer $n$
such that $\Delta (P_\infty )\cap F_n$ is not meager in $\Delta (P_\infty )$. Therefore there 
exists a sequence $s$ in ${\cal S}\setminus\{\emptyset\}$ such that ${\Delta
(N_s\cap P_\infty)\subseteq F_n}$. Thus $\Delta (N_s)\subseteq F_n$ and $(s^\frown
0^\omega,s^\frown 0^\omega)\in
\phi^{-1}(\{ 0\})$, which is absurd.\bigskip

\noindent (b) Let $A:= I[P_\infty]$. As $I$ is a homeomorphism from $2^\omega$ onto its range and
$P_\infty\approx\omega^\omega$, we have $A\approx\omega^\omega$. We have $F := \phi^{-1}(\{ 1\})\cap 
(2^\omega\times A)\equiv\bigcup_{s\in {\cal S}}~\{s^\frown
0^\omega\}\times (N_s\cap P_\infty )$. Let us show that $\overline{F}^{2^\omega\times
P_\infty}\subseteq F\cup\Delta (P_\infty )$. Then $F =
\overline{F}^{2^\omega\times P_\infty}\setminus\Delta (P_\infty )$ will be 
$D_2(\boraone)(2^\omega\times A)\subseteq \bortwo(2^\omega\times A)$. As $\phi^{-1}(\{
0\})\cap (2^\omega\times A) = (2^\omega\times A)\setminus \phi^{-1}(\{ 1\})$, we will
have $\phi\in {\cal B}_1(2^\omega\times A,2)$. If $(s_n^\frown 0^\omega,s_n^\frown\gamma_n)\in
F$ tends to $(\beta,\alpha)\in
(2^\omega\times P_\infty )\setminus F$, we may assume that $\vert s_n\vert $ increases strictly. So
for each integer $p$ and for $n$ big enough we have $\beta
(p) = s_n(p) =\alpha (p)$. Thus $\alpha =\beta$.\bigskip

 If $A$ were uniformly recoverable, we could find a dense sequence 
 $D:=(x_p)$ of 
$2^\omega$ such that every function of $A$ is recoverable with respect
to $(x_p)$. Let $s\in {\cal S}$. Then $I(s^\frown 1^\omega)$ is in 
$A$, and it is the characteristic function of the following set:
$${\{s\lceil n^\frown 0^\omega~/~n=0~~\mbox{or}~~(0<n\leq\vert s\vert ~~\mbox{and}~~s(n-1)=1)\}\cup
\{s^\frown 1^{p+1}0^\omega~/~p\in\omega\}}.$$
For $n$ big enough,
$s_n[s^\frown 0^\omega,D]$ is in this set, thus $s^\frown 0^\omega\in
D$ and $P_f:=2^\omega\setminus
P_\infty\subseteq D$. So the functions of
$I[2^\omega]$ are all recoverable with respect to $D$. But this
contradicts the previous point.$\hfill\square$\bigskip

\noindent\bf (B) Study~of~the~link~between~recoverability~and~ranks~on~Baire~class~one~functions.\rm\bigskip

  So there exists a metrizable compact set of characteristic functions of $D_2(\boraone)$
sets which is not uniformly recoverable. So the boundedness of the complexity of
functions of
$A$ does not insure that $A$ is uniformly recoverable. Notice that the example of the ``split
interval" is another proof of this, in the case where the compact space is not metrizable.
Indeed, functions of the ``split interval" are characteristic functions of open or closed
subsets of
$[0,1]$ (of the form $]a,1]$ or $[a,1]$, with $a\in [0,1]$).\bigskip

 In [B2], the author
introduces a rank which measures the complexity of numeric Baire class one functions defined
on a metrizable compact space. Let us recall this definition, which makes sense for functions
defined on a Polish space $X$ which is not necessarily compact.

\vfill\eject

\noindent $\bullet$ Let ${\cal A}$ and ${\cal B}$ be two disjoint $G_\delta$ subsets of $X$, and $R({\cal
A},{\cal B})$ be the set of increasing sequences $(G_\alpha
)_{\alpha\leq\beta}$ of open subsets of $X$, with $\beta<\omega_1$, which satisfy\bigskip 

\noindent 1. $G_{\alpha +1}\setminus G_\alpha$ is disjoint from ${\cal A}$ or from ${\cal B}$ if $\alpha
<\beta$.\smallskip

\noindent 2. $G_\gamma = \cup_{\alpha <\gamma}~G_\alpha$ if $0<\gamma\leq\beta$ is a limit ordinal.\smallskip
 
\noindent 3. $G_0=\emptyset$ and $G_\beta = X$.\bigskip
 
 Then $R({\cal A},{\cal B})$ is not empty, because ${\cal A}$ and ${\cal B}$ can be separated
by a $\bortwo$ set, which is of the form 
$$D_\xi\big((U_\alpha)_{\alpha<\xi}\big) :=
\bigcup_{\alpha<\xi~\mbox{with~parity~opposite~to~that~of}~\xi} 
U_\alpha\setminus\big(\cup_{\theta <\alpha} U_\theta\big)\mbox{,}$$ 
where $(U_\alpha)_{\alpha<\xi}$ is an increasing sequence of open subsets of $X$ and 
$1\leq\xi<\omega_1$ (see [Ke]). Then we check that $(G_\alpha
)_{\alpha\leq\xi +1}\in R({\cal A},{\cal B})$, where 
$G_{\alpha +1} := U_\alpha~\mbox{if}~\alpha <\xi$.\bigskip

\noindent $\bullet$ We set $L({\cal A},{\cal B}) := \mbox{min}\{\beta <\omega_1~/~\exists~(G_\alpha
)_{\alpha\leq\beta}\in R({\cal A},{\cal B})\}$. If $f\in {\cal B}_1(X,\mathbb{R})$ and $a<b$ are
real numbers, we let  
$L(f,a,b) := L(\{f\leq a\},\{f\geq b\})$. Finally, 
$$L(f) := \mbox{sup}\{L(f,q_1,q_2)~/~q_1<q_2\in\mathbb{Q}~\}.$$
In [B2], the author shows that, if $A\subseteq {\cal C}(X,\mathbb{R})$ is relatively compact in 
${\cal B}_1(X,\mathbb{R})$, then
$${\mbox{sup}\{L(f,a,b)~/~f\in\overline{A}^{\mbox{p.c.}}\}<\omega_1}$$ 
if $X$ is a compact space and if $a<b$ are real numbers. He wonders whether his result
remains true for a separable compact subspace $A$ of ${\cal B}_1(X,\mathbb{R} )$.\bigskip

 We can ask the question of the link between uniform recoverability
of $A$ and the fact that 
$$\mbox{sup}\{L(f)/f\in A\}\! <\! \omega_1.$$ 
If ${D_\xi(\boraone)(X) :=\{D_\xi\big((U_\eta)_{\eta<\xi}\big)/(U_\eta)_{\eta<\xi}\subseteq\boraone(X)~
\mbox{increasing}\}}$ and $A\in D_\xi(\boraone)(X)$, one has $\check {\cal A}\in D_{\xi+1}(\boraone)(X)$ and  $L(\check {\cal A},{\cal A})\leq \xi +2$ by the previous facts. So the rank of 
the characteristic function of ${\cal A}$ is at most $\xi +2$. In the case of the example in
Corollary 16 and of the ``split interval", one has  
$\mbox{sup}\{L(f)~/~f\in A\}\leq 4<\omega_1$. Therefore, the fact that $L$ is bounded on $A$
does not imply uniform recoverability of $A$, does not imply that $\phi$
is Baire class one, and does not imply that $A$ is metrizable. But we have the following
result. It is a partial answer to J. Bourgain's question.

\begin{prop} If $X$ is a Polish space,
$Y\subseteq\mathbb{R}$ and $A\subseteq {\cal B}_1(X,Y)$ is a Polish space, then we have 
$\mbox{sup}\{L(f)~/~f\in A\}<\omega_1$.\end{prop} 

\noindent\bf Proof.\rm\ Let $a<b$ be real numbers, ${\cal A} :=
\{(x,f)\in X\times A~/~f(x)\leq a\}$ and 
$${{\cal B} :=\{(x,f)\in X\times A~/~f(x)\geq b\}}.$$
As $\phi$ is Baire class two, ${\cal A}$ and ${\cal B}$ are $\bormthree (X\times A)$ with horizontal sections in $\bormtwo(X)$.

\vfill\eject

 So there exists a finer Polish topology
$\tau_{\cal A}$ on $A$ such that 
${\cal A}\in\bormtwo(X\times [A,\tau_{\cal A}])$ (see [L1]). The same thing is true
for ${\cal B}$. Let
$\tau$ be a Polish topology on $A$, finer than $\tau_{\cal A}$ and $\tau_{\cal B}$ (see 
Lemma 13.3 in [Ke]). As ${\cal A}$ and ${\cal B}$ are disjoint, there exists
$\Delta_{a,b}\in\bortwo(X\times [A,\tau])$ which separates ${\cal A}$ from ${\cal B}$. Let
$\xi_{a,b} <\omega_1$ be such that $\Delta_{a,b}\in D_{\xi_{a,b}}(\boraone)(X\times
[A,\tau])$. For each function $f$ of $A$, the set 
$\Delta_{a,b}^f$ is a $D_{\xi_{a,b}}(\boraone)(X)$ which separates $\{f\leq a\}$ from $\{f\geq
b\}$. Thus $L(\{f\leq a\},\{f\geq b\})\leq \xi_{a,b} +1$. Therefore 
$\mbox{sup}\{L(f)~/~f\in A\}\leq\mbox{sup}\{L(f,a,b)~/~a<b\in\mathbb{Q} ~\}<\omega_1$.$\hfill\square$

\begin{cor} If $X$ is a Polish space,
$Y\subseteq\mathbb{R}$ and if $A\subseteq {\cal B}_1(X,Y)$ is uniformly recoverable and
compact, then $\mbox{sup}\{L(f)~/~f\in A\}<\omega_1$.\end{cor}

 We can wonder whether this result is true for the set of recoverable
functions with respect to a dense sequence of $X$. We will see that it is not the case.

\begin{prop} Let $(x_p)$ be a dense sequence of a nonempty 
perfect Polish space $X$, and $Y := 2$. Then $\mbox{sup}\{L(f)~/~f~
\mbox{is~recoverable~with~respect~to}~(x_p)\}=\omega_1$.\end{prop}

\noindent\bf Proof.\rm\ Set $D$ of the elements of the dense sequence is
countable, metrizable, nonempty and perfect. Indeed, if $x_p$ is an isolated point of
$D$, then it is also isolated in $X$, which is absurd. Thus $D$ is homeomorphic to 
$\mathbb{Q} ~$ (see 7.12 in [Ke]). For $1\leq\xi <\omega_1$, there exists a countable
metrizable compact space $K_\xi$ and ${\cal A}_\xi\in D_\xi (\boraone)(K_\xi )\setminus\check
D_\xi (\boraone)(K_\xi )$ (see [LSR]). So we may assume that $K_\xi\subseteq D$ (see
7.12 in [Ke]). Thus we have ${\cal A}_\xi\notin \check D_\xi (\boraone)(X)$. We will
deduce from this the fact that $L(1\!\! {\rm I}_{{\cal A}_{\xi +1}})>\xi$.\bigskip

 To see this, let us show that, if $L(1\!\! {\rm I}_{\cal A})=L(\check {\cal
A},{\cal A})\leq\xi$, then ${\cal A}\in
\check D_{\xi +1}(\boraone)(X)$. Let ${(G_\alpha)_{\alpha\leq\xi'}}$ be in ${R(\check {\cal
A},{\cal A})}$, where 
$\xi'\in\{\xi,\xi+1\}$ is odd. We let, for $\alpha <\xi'$,
$$U_\alpha :=
\left\{\!\!\!\!\!\!\!\!
\begin{array}{ll} 
& \bigcup_{\theta <\alpha} U_\theta\cup\bigcup_{\theta\leq\alpha/{\cal A}\cap G_{\theta +1}\setminus G_\theta =
\emptyset}  G_{\theta +1}~~\mbox{if}~\alpha~\mbox{is~even,}\cr\cr 
& \bigcup_{\theta <\alpha} U_\theta\cup\bigcup_{\theta\leq\alpha /\check {\cal A}\cap G_{\theta +1}\setminus G_\theta =
\emptyset}  G_{\theta +1}~~\mbox{if}~\alpha~\mbox{is~odd.}
\end{array}
\right.$$
Then $D_{\xi'}\big((U_\alpha)_{\alpha <\xi'}\big)$ separates $\check {\cal A}$ from ${\cal A}$.
Indeed, if $x\notin {\cal A}$, let $\alpha'\leq\xi'$ be minimal such that $x\in G_{\alpha'}$.
Then $\alpha'$ is the successor of $\alpha<\xi'$, and $x\in \check {\cal A}\cap G_{\alpha
+1}\setminus G_\alpha$. So ${\cal A}\cap G_{\alpha +1}\setminus G_\alpha = \emptyset$, by 
condition 1. If
$\alpha$ is even, then $x\in U_\alpha\setminus (\cup_{\theta <\alpha} U_\theta )$ because $U_\theta\subseteq
G_{\theta +1}$ if $\theta <\xi'$. If $\alpha$ is odd, then $x\in U_{\alpha +1}\setminus U_\alpha$.
In both cases, $x\in D_{\xi'}\big((U_\alpha)_{\alpha <\xi'}\big)$. If $x\in
U_\alpha\setminus (\cup_{\theta <\alpha} U_\theta )$ with $\alpha <\xi'$ even, there exists
$\theta\leq\alpha$ such that $x\in G_{\theta +1}$ and ${\cal A}\cap G_{\theta +1}\setminus G_\theta =
\emptyset$. Let $\eta'\leq\xi'$ be minimal such that 
$x\in G_{\eta'}$. As before, $\eta'$ is the successor of $\eta <\xi'$. Let us argue 
by contradiction: we assume that $x\in {\cal A}$. Then $x\in {\cal
A}\cap G_{\eta +1}\setminus G_\eta\not=
\emptyset$, so $\check {\cal A}\cap G_{\eta +1}\setminus G_\eta =\emptyset$. If $\eta$ is odd,
then $x\in U_\eta$, thus $\eta = \alpha$. This contradicts the parity of $\alpha$. If 
$\eta$ is even, then $x\in U_{\eta +1}$ and $\eta =\alpha = \theta$. So $x\in G_{\theta +1}\setminus G_\theta \subseteq \check {\cal A}$. This is the contradiction we were looking
for.\bigskip

 It remains to check that $1\!\! {\rm I}_{{\cal A}_{\xi +1}}$ is recoverable. If
$x\in D$, then ${s_n[x,D]=x}$ for almost all integer $n$. Thus $1\!\! {\rm
I}_{{\cal A}_{\xi +1}}\big(s_n[x,D]\big)$ tends to $1\!\! {\rm I}_{{\cal A}_{\xi
+1}}(x)$. If $x$ is not in $D$, then ${x\notin {\cal A}_{\xi +1}\subseteq K_{\xi
+1}\subseteq D}$. So, from some point on, $s_n[x,D]\notin K_{\xi +1}$, and $1\!\! {\rm
I}_{{\cal A}_{\xi +1}}\big(s_n[x,D]\big)$ is ultimately constant and tends to $1\!\! {\rm
I}_{{\cal A}_{\xi +1}}(x)$.$\hfill\square$\bigskip

\noindent\bf Remark.\rm ~We can find in [KL] the study of some other
ranks on Baire class one functions. The rank $L$ is essentially the separation rank defined
in this paper. In the case where $X$ is a metrizable compact space and where the Baire class
one functions considered are bounded, Propositions 17, 19 and Corollary 
18 remain valid for these other ranks.

\vfill\eject

\noindent\bf (C) Sufficient~conditions~for~uniform~recoverability.\rm

\begin{thm} Assume that $Y$ is a metric space, and that $A$, equipped
with the compact open topology, is a separable subset of ${\cal B}_1(X,Y)$. Then $A$ is
uniformly recoverable.\end{thm}

\noindent\bf Proof.\rm\ Let $(l_q)$ be a dense sequence of $A$ for the compact
open topology. By the lemma showed in [Ku], page 388, for each integer $q$ there exists a sequence
$(h_n^q)_n\subseteq {\cal B}_1(X,Y)$  which uniformly tends to $l_q$, functions $h_n^q$ having
a discrete range. Enumerating the sequence $(h_n^q)_{n,q}$, we get $(h_n)_{n}$. Every function of $A$ is in the closure of this sequence for the compact open topology. For each integer $n$, one can get a 
countable partition $(B^n_{p})_{p}$ of $X$ into $\bortwo$ sets on which $h_{n}$ is constant. 
Express each of these sets as a countable union of closed sets. Putting 
all these closed sets together gives 
a countable sequence of closed subsets of $X$. As in the proof of Theorem 4, this gives $D$ which 
approximates each of these closed sets. Now let $f\in A$, $x\in X$ and $\varepsilon >0$. Consider the 
compact subset $K := {\cal R}\big(x,D\big)\cup\{ x\}$ of $X$. By uniform convergence on $K$, there 
exists $N\in\omega$ such that, for each $t$ in $K$, we have ${d_{Y}(f(t),h_N(t))<\varepsilon /2}$. Let $p$ be 
an integer such that $x\in B^N_{p}$. Now $K\setminus B^N_{p}$ is finite and we have 
$h_{N}(s_{n}[x,D])=h_{N}(x)$ for each $n\in\omega$, except maybe a finite 
number of them. So we have the following inequality, for all but finitely 
many $n$:
$$\begin{array}{ll} 
& d_{Y}\big(f(x),f(s_{n}[x,D])\big) \leq\cr 
& ~~~~~~~~~~~~d_{Y}\big(f(x),h_N(x)\big)+d_{Y}\big(h_N(x),h_N(s_{n}[x,D])\big)+d_{Y}\big(h_N(s_{n}[x,D]),f(s_{n}[x,D])\big)\! <\! \varepsilon 
\end{array}$$ 
(this last argument is essentially in [DE]).$\hfill\square$\bigskip 

  The following corollary has been showed in [FV] when $X=\mathbb{R}$ 
 and with another way of extracting the subsequence.

\begin{cor} Let $A\subseteq {\cal B}_1(X,Y)$ be countable. Then 
$A$ is uniformly recoverable.\end{cor}

\noindent\bf Proof.\rm\ Put a compatible distance on 
$Y$.$\hfill\square$

\begin{prop} Let $(Y_p)$ be a basis for the topology of
$Y$, and\medskip

\noindent (1) For each integer $p$, $\phi^{-1}(Y_p) \in \big(\bormone (X)\times {\cal
P}(A)\big)_\sigma$.\smallskip

\noindent (2) There exists a finer metrizable separable topology on $X$, made of $\boratwo(X)$, and
making functions of $A$ continuous.\smallskip

\noindent (3) $A$ is uniformly recoverable.\medskip

 Then (1) $\Leftrightarrow$ (2) $\Rightarrow$ (3).\end{prop}

\noindent\bf Proof.\rm\ (1) $\Rightarrow$ (2) We have $\phi^{-1}(Y_p) =
\bigcup_{n\in\omega} F^p_n\times B^p_n$, where $F^p_n$ is a closed subset of $X$ and
${B^p_n\subseteq A}$. If $f\in A$, then $f^{-1}(Y_p) = \phi^{-1}(Y_p)^f =
\bigcup_{n/f\in B^p_n} F^p_n$. Therefore, it is enough to find a finer metrizable separable
topology on 
$X$, made of $\boratwo(X)$, and making the $F^p_n$'s open. Let $(X_n)$ be a 
basis for the topology of $X$, closed under finite intersections, and $(G_q)$ be the sequence
of finite intersections of $F^p_n$'s. Then set $\tau$ of unions of sets of the form $X_n$ or
$X_n\cap G_q$ is a topology, with a countable basis, made of $\boratwo(X)$, finer than the
initial topology on $X$ (thus Hausdorff), and makes the $F^p_n$'s open. It remains to check
that it is regular.

\vfill\eject

 So let $x\in X$ and $F \in\bormone (X,\tau)$, with $x\notin F$. We have 
$X\setminus F =\bigcup_{p} X_{n_p} \cup \bigcup_{k} X_{m_k}\cap G_{q_k}$. Either there exists $p$ such that $x\in X_{n_p}$; in this case, by regularity of initial topology on $X$ we can find
two disjoint open sets $V_1$ and $V_2$ with $x\in V_1$ and $X\setminus X_{n_p}\subseteq V_2$.
But these two open sets are $\tau$-open and $F\subseteq X\setminus
X_{n_p}\subseteq V_2$. Or there exists $k$ such that $x\in X_{m_k}\cap G_{q_k}$; in this case,
by regularity of initial topology on $X$, we can find two disjoint open sets $W_1$ and
$W_2$ with $x\in W_1$ and $X\setminus X_{m_k}\subseteq W_2$. But then $W_1\cap G_{q_k}$ and
$W_2\cup (X\setminus G_{q_k})$ are
$\tau$-open and disjoint, $x\in W_1\cap G_{q_k}$ and $F\subseteq (X\setminus X_{m_k}) \cup
(X\setminus G_{q_k}) \subseteq W_2 \cup (X\setminus G_{q_k})$.\bigskip

\noindent (2) $\Rightarrow$ (3) Let $\tau$ be the finer topology. Then identity map from $X$,
equipped with its initial topology, into $X$, equipped with $\tau$, is Baire class one.
Therefore, it is recoverable. So let $(x_p)$ be a
dense sequence of $X$ such that, for each $x\in X$, $s_n[x,(x_p)]$ tends to $x$, in the sense of
$\tau$. Let
$f\in A$. As $f$ is continuous if $X$ is equipped with $\tau$, $f\big(s_n[x,(x_p)]\big)$ tends
to $f(x)$ for each $x\in X$. Therefore $f$ is recoverable with respect to
$(x_p)$.\bigskip

\noindent (2) $\Rightarrow$ (1) Let $(X_n)$ be a basis for finer topology $\tau$ (therefore, we have
$X_n\in\boratwo(X)$). Let 
$$C^p_n := \{f\in A~/~X_n\subseteq\phi^{-1}(Y_p)^f\}.$$ 
Then $\phi^{-1}(Y_p) =\bigcup_{n} X_{n}\times C^p_{n} \in \big(\bormone (X)\times {\cal
P}(A)\big)_\sigma$.$\hfill\square$\bigskip 

\noindent\bf Remark.\rm ~If $X$ is a standard Borel space and $A$ is a Polish space, conditions
(1) and (2) of Proposition 22 are equivalent to ``For each integer $p$, $\phi^{-1}(Y_p) \in \big(\bormone (X)
\times \borel (A)\big)_\sigma$". Indeed, let $P$ be a Polish space such that $X$ is a Borel
subset of $P$, and $f\in A$. As $f$ is continuous if $X$ is equipped with $\tau$, $f^{-1}(Y_p)
=\bigcup_{k} X_{n^{p,f}_k}$ for each integer $p$. Let $C^p_n := \{f\in A~/~X_n\subseteq
\phi^{-1}(Y_p)^f\}$. Then $C^p_n$ is $\ca (A)$, because $\phi$ is Baire class two: 
$$f\in C^p_n\ \Leftrightarrow\ \forall x\in P~~x\notin X_n~~\mbox{or}~~\phi (x,f)\in Y_p.$$
Moreover, $\phi^{-1}(Y_p) = \bigcup_{n} X_{n}\times C^p_{n}$. By $\borel$-selection
(see 4B5 in [M]), there exists a Borel function $N_p:P\times
A\rightarrow\omega$ such that $(x,f)\in X^{}_{N_p(x,f)}\times
C^p_{N_p(x,f)}$ if $f(x)\in Y_p$. Let 
$$S^p_{n} := \{f\in A~/~\exists x\in X~~N_p(x,f)=n~~\mbox{and}~~\phi(x,f)\in Y_p\}.$$ 
Then $S^p_{n}\in \ana (A)$ and is a subset of $C^p_{n}$ ; by the separation therem, there exists a Borel subset $B^p_{n}$ of $A$ such that $S^p_{n}\subseteq B^p_{n}\subseteq C^p_{n}$. Then we have 
$\phi^{-1}(Y_p) =\bigcup_{n} X_{n}\times B^p_{n} \in \big(\bormone (X)\times \borel (A)\big)_\sigma$.

\begin{prop} If $A$ has a countable basis, then there exists a finer
metrizable separable topology on $X$ making the functions of $A$ continuous. Moreover, if $X$ is
Polish, we can have this topology Polish.\end{prop}

\noindent\bf Proof.\rm\ Let $(A_n)$ be a basis for the topology of
$A$, and
${X_n^p :=\{x\!\in\! X/A_n\subseteq\phi^{-1}\! (Y_p)_x\}}$. As 
$$\phi^{-1}(Y_p)_x =\{f\in A~/~f(x)\in Y_p\}\in\boraone(A)\mbox{,}$$ 
we have $\phi^{-1}(Y_p) = \bigcup_{n\in\omega} X_n^p\times A_n$ and $f^{-1}(Y_p) = \phi^{-1}(Y_p)^f = \bigcup_{n/f\in A_n} X_n^p$ for each
$f\in A$. Thus it is enough to find a finer metrizable separable topology on $X$ making 
$X_n^p$'s open.

\vfill\eject

 We use the same method as the one used to prove implication
(1) $\Rightarrow$ (2) of Proposition 22. We notice that the algebra generated by the $X_n^p$'s is
countable (we let $(G_q)$ be the elements of this algebra).\bigskip

 As $\phi$ is Baire class two, $\phi^{-1}(Y_p)$ is a 
$\borathree $ set with vertical sections in $\boraone(A)$. If $X$ and $A$ are Polish, we deduce from [L1] the existence of a finer Polish topology $\tau_p$ on $X$ such that 
$${\phi^{-1}(Y_p)\in\big(\boraone (X,\tau_p)\times \boraone (A)\big)_\sigma}.$$ 
Let $(B^p_n)_n$ be a basis for $\tau_p$. Then there  exists a finer Polish topology $\tau$ on $X$ making the Borel sets $B^p_n$'s open (see  Exercises 15.4 and 13.5 in [Ke]). Then we are done, because $\tau$
is finer than the $\tau_p$'s.$\hfill\square$\bigskip

 Therefore, the problem is to find the finer topology in $\boratwo(X)$. We have seen that
it is not the case in general. If we look at
Propositions 15 and 22, we can wonder whether conditions of Proposition 
22 and the fact that
$\phi$ is Baire class one are equivalent, especially in the case where $Y$ is 0-dimensional.
This question leads to the study of Borel subsets of $2^\omega\times2^\omega$. The answer is no in
general. First, because of Corollary 16. It shows that the fact that $\phi$ is Baire class one
does not imply uniform recoverability (with $A$ Polish, in fact homeomorphic to $\omega^\omega$). 
Secondly, let $A :=
\{f\in {\cal
B}_1(2^\omega,2)~/~f~\mbox{is~recoverable~with~respect~to}~(x_p)\}$, where 
$(x_p):=P_f$ is dense in $2^\omega$. Then $A$ is uniformly recoverable, but we cannot find a finer 
metrizable separable topology $\tau$ on $2^\omega$, made of $\boratwo(2^\omega)$ and making 
the functions of $A$ continuous. Otherwise, the characteristic functions of the compact sets ${K_x:=\{
x\}\cup
\{s_n[x,(x_p)]~/~n\in\omega\}}$ would be continuous for $\tau$, and this would contradict
the Lindel\"of property, with $\bigcup_{x\in P_\infty} K_x$. But $A$ has no
countable basis. Otherwise, the set of charateristic functions of the sets 
$K_x$ (for $x\in P_\infty$) would also have one; this would contradict the Lindel\"of property
too (this last set is a subset of $\bigcup_{x\in P_\infty}
\{f\in {\cal B}_1(2^\omega ,2)~/~f(x)=1\}$). This leads us to assume that $A$ is a $K_\sigma$
and metrizable space, to hope for such an equivalence.\bigskip  

 If $\phi$ is Baire class one, then $\phi^{-1}(Y_p)$ is a $\boratwo$ subset of $X\times A$ with 
vertical sections in $\boraone(A)$. Thus it is natural to ask the\bigskip

 \noindent\bf Question.\rm ~Does every $\boratwo$ subset of $X\times A$ with vertical sections
in $\boraone(A)$ belong to the class $\big(\bormone (X)\times {\cal P}(A)\big)_\sigma$?\bigskip

 If the answer is yes, then the fact that $\phi$ is Baire class one implies condition (1) in
Proposition 22, and the conditions of this proposition are equivalent to the fact that $\phi$ is
Baire class one. The answer is negative, even if we assume that $X$ and $A$ are metrizable
compact spaces:

\begin{prop}  There exists a $\check D_2(\boraone)$ subset of 
$2^\omega\times2^\omega$ with vertical sections in $\borone (2^\omega)$ which is not $\big(\bormone 
(2^\omega)\times {\cal P}(2^\omega)\big)_\sigma$.\end{prop}

\noindent\bf Proof.\rm\ Let $E :=
(P_\infty\times2^\omega )\cup\bigcup_{s\in {\cal S}}~\{s^\frown
0^\omega\}\times (\check N_s\cup N_{s^\frown 0})$ (we use again notations of the proof of
Corollary 16). Clearly, vertical sections of $E$ are $\borone (2^\omega)$. We set 
$$G:=\{\alpha\in2^\omega~/~\forall n~\exists m\geq n~~\alpha (m) = \alpha
(m+1) = 1\}.$$ 
This is a dense $G_\delta$ subset of $2^\omega$, included in $P_\infty$. If 
$\alpha\notin G$, then the horizontal section $\check E^\alpha$ is finite.

\vfill\eject

 Otherwise, it is infinite and 
countable (it is a subset of $P_f$), and it is a sequence which tends to $\alpha$. If 
$(s_n^\frown 0^\omega,s_n^\frown\gamma_n)_n\subseteq\check E$ tends to $(\beta,\alpha)$, then there are essentially two cases. Either the length of $s_n$ is strictly increasing and $\alpha =\beta$.
Or we may assume that $(s_n)$ is constant and $(\beta,\alpha)\notin E$. As diagonal 
${\Delta (2^\omega )\subseteq E}$, we can deduce from this that ${\check E =\overline{\check
E}\setminus
\Delta (2^\omega )\in D_2(\boraone)(2^\omega\times2^\omega)}$. Assume that $E\in
\big(\bormone (2^\omega)\times {\cal P}(2^\omega)\big)_\sigma$. We have $E = \bigcup_n F_n\times
E_n$, where $F_n\in\bormone (2^\omega)$ and $E_n\subseteq2^\omega$. Let $C_n :=
\{\alpha\in2^\omega~/~F_n\subseteq E^\alpha\}$. Then $C_n\in\ca (2^\omega )$ and $E = \bigcup_n
F_n\times C_n$. As
$\Delta (2^\omega )\subseteq E$, $2^\omega\subseteq \bigcup_n F_n\cap C_n$. So there exists an
integer
$n$ such that $F_n\cap C_n$ is not meager, and a sequence $s\in 2^{<\omega}$ such that 
$N_s\cap F_n\cap C_n$ is a comeager subset of
$N_s$. In particular,
$N_s\subseteq F_n$. As $G$ is comeager, there exists $\alpha\in G\cap N_s\cap C_n$.
Let $(\beta_m)\subseteq \check E^\alpha$ converging to $\alpha$. From some point $m_0$ on, we
have $\beta_m\in N_s$. So $(\beta_m,\alpha)\in F_n\times C_n\subseteq E$ if
$m\geq m_0$. But this is absurd because $(\beta_m,\alpha)\notin E$.$\hfill\square$\bigskip

 We can specify this result:
 
\begin{prop} There exists a metrizable compact space 
$A\subseteq {\cal B}_1(2^\omega ,2)$ which is uniformly recoverable, but for which we cannot
find any finer metrizable separable topology on $2^\omega$, made of $\boratwo (2^\omega )$, making the 
functions of $A$ continuous.\end{prop}

\noindent\bf Proof.\rm\ We use again the notation of the proof 
of Proposition 24. Let
$\psi:\omega\rightarrow {\cal S}$ be a bijective map such that for $s,t\in
{\cal S}$,
$s\prec_{\not=} t~$  implies $\psi^{-1}(s)<\psi^{-1}(t)$. Such a bijection exists. Indeed, we 
take $\psi := (\theta ~\circ ~\phi_{\lceil {\cal S}})^{-1}$, where 
$\theta :\phi[{\cal S}]\rightarrow \omega$ is an increasing bijection, and where 
$$\phi :\left\{\!\!
\begin{array}{ll}
2^{<\omega} 
& \!\!\!\!\rightarrow \omega\cr 
s 
& \!\!\!\!\mapsto\left\{\!\!\!\!\!\!\!\!
\begin{array}{ll}
& 0~~\mbox{if}~~s=\emptyset\mbox{,}\cr\cr 
& q_0^{s(0)+1}...q_{s(\vert s\vert -1)}^{s(\vert s\vert -1)+1}~\mbox{otherwise.}
\end{array}
\right.
\end{array}
\right.$$ 
(where $(q_n)$ is sequence of prime numbers). We let $x_{2n}:=\psi
(n)^\frown 1^\omega$, $x_{2n+1}:=\psi (n)^\frown 0^\omega$, and $x_s:=x_{\mbox{min}\{p\in\omega/s\prec x_p\}}$ if $s\in 2^{<\omega}$.\bigskip

\noindent $\bullet$ Let us show that, if $s\in 2^{<\omega}\setminus\{\emptyset\}$, then $x_s\in
P_\infty$ is equivalent to $s\in {\cal S}$. If $s\in {\cal S}$ and $x_s\in P_f$, then there
exists $u$ in $\cal S$ such that $x_s = s^\frown u^\frown 0^\omega$. Then 
$x_{2\psi^{-1}(s)}$ comes strictly before $x_{2\psi^{-1}(s)+1}$, which comes
before $x_{2\psi^{-1}(s^\frown u)+1} = x_s$. But $s\prec s^\frown 1^\omega = x_{2\psi^{-1}(s)}$,
which is absurd.\bigskip 

 If $s\notin {\cal S}$ and $x_s\in P_{\infty}$, there exists $u$ in $\cal S$ such that $x_s =
s^\frown (1^{\vert u\vert }-u)^\frown 1^\omega$. Let
$s'\in {\cal S}$ and $m$ be an integer such that $s = s'^\frown 0^{m+1}$. Then 
$x_{2\psi^{-1}(s')+1}$ comes strictly before $x_{2\psi^{-1}(s^\frown
(1^{\vert u\vert }-u)^\frown 1)}$, which comes before $x_s$. But $s\prec s^\frown 0^\omega =
s'^\frown 0^\omega =x_{2\psi^{-1}(s')+1}$, which is absurd.\bigskip

\noindent $\bullet$ We set 
$$I\! :\! \left\{\!\!
\begin{array}{ll}
2^\omega 
& \!\rightarrow\!  {\cal B}_1(2^\omega, 2)\cr\cr 
\alpha 
& \!\mapsto\!  \left\{
\begin{array}{ll}
2^\omega 
& \!\!\!\!\rightarrow  2\cr 
\beta 
& \!\!\!\!\mapsto\left\{\!\!\!\!\!\!\!\!
\begin{array}{ll}
& 0~~\mbox{if}~~\exists s\in {\cal S}~~\beta\!  =\!  s^\frown 0^\omega~~\mbox{and}~~\alpha\!
\in\!  N_{s^\frown 1}\mbox{,}\cr\cr 
& 1~~\mbox{otherwise.}
\end{array}
\right.
\end{array}
\right.
\end{array}
\right.$$

\vfill\eject

 Then $I$ is defined because $\{\beta\in2^\omega~/~\beta\notin I(\alpha )\}$ is $\overline{\{
 \beta\in2^\omega~/~\beta\notin I(\alpha )\}}\setminus\{\alpha\}\in D_2(\boraone)(2^\omega )$ if
$\alpha\in G$, and is finite otherwise. $I$ is continuous because 
$$I(\alpha )(\beta)=1\Leftrightarrow \left\{\!\!\!\!\!\!\!\!
\begin{array}{ll} 
& \alpha\in2^\omega~~\mbox{if}~~\beta\in P_\infty\mbox{,}\cr\cr 
& \check N_s\cup N_{s^\frown 0}~~\mbox{if}~~\beta=s^\frown 0^\omega~~\mbox{and}~~s\in {\cal S}.
\end{array}
\right.$$
 Therefore, $A:=I[2^\omega ]$ is an analytic compact space and is metrizable.\bigskip
 
  As $E = (\mbox{Id}_{2^\omega}\times I)^{-1}\big(\phi^{-1}(\{ 1\})\big)$, $\phi^{-1}(\{ 1\})\notin \big(\bormone (2^\omega )\times {\cal P}(A)\big)_\sigma$, by Proposition 
24. So there
is no finer metrizable separable topology on $2^\omega$, made of 
$\boratwo(2^\omega )$ and making the functions of $A$ continuous, by Proposition 
22. But 
$A$ is uniformly recoverable with respect to $(x_p)$.\bigskip

 Indeed, as $P_f\subseteq (x_p)$, it is enough to see that if $\alpha\in G$, then 
$I(\alpha)$ is recoverable with respect to $D:=(x_p)$. The only thing to see is that
from some integer $n_0$ on, ${s_{n}[\alpha,D]\in E^\alpha}$. We may assume that
$\alpha\notin D$ because $G\subseteq P_\infty$.\bigskip

 We take $(W_m) := (N_s)_{s\in
2^{<\omega}}$ as a good basis for the topology of $2^\omega$. So that, if $\alpha\notin D$, 
$$\begin{array}{ll}
s_{n+1}[\alpha, D] 
& = x_{\mbox{min} \big\{ p\in\omega/\exists s\in 2^{<\omega}~\alpha, x_p\in N_s\subseteq2^\omega\setminus\{s_{0}[\alpha,D],...,s_{n}[\alpha,D]\} \big\} } \cr 
& = x_{\mbox{min} \big\{ p\in\omega/\alpha\lceil (\mbox{max}_{q\leq n} \vert\alpha\wedge
s_{q}[\alpha,D]\vert +1)\prec x_p \big\} }.
\end{array}$$
But as the sequence $(\vert\alpha\wedge s_{n}[\alpha,D]\vert )_n$ is strictly increasing, 
$\mbox{max}_{q\leq n} \vert\alpha\wedge s_{q}[\alpha,D]\vert =\vert\alpha\wedge s_{n}[\alpha,D]\vert$.  
Thus $s_{n+1}[\alpha, D]=x_{\alpha\lceil (\vert\alpha \wedge s_{n}[\alpha, D]\vert +1)}$. By the previous facts, it is enough to get 
$${x_{\alpha\lceil (\vert\alpha\wedge s_{n}[\alpha, D]\vert +1)}\! \in\! E^\alpha}.$$ 
Let $M_n\!  :=\!  \vert\alpha\wedge\!  s_{n}[\alpha, D]\vert$. If $\alpha (M_n)=1$, then ${s_{n+1}[\alpha, D]}$ is in 
${P_\infty\subseteq E^\alpha}$. Otherwise, 
$${s_{n+1}[\alpha, D]=\alpha\lceil(M_n+1)^\frown u^\frown 0^\omega}\mbox{,}$$ 
where $u\in {\cal S}$. If $u\not=\emptyset$, then $s_{n+1}[\alpha, D]$ is minimal in 
$N_{\alpha\lceil(M_n+1)^\frown u}\subseteq N_{\alpha\lceil(M_n+1)}$, so 
$s_{n+1}[\alpha,D]$ is in $P_\infty$, which is absurd. Thus $u=\emptyset$ and 
$s_{n+1}[\alpha,D]\in E^\alpha$.$\hfill\square$\bigskip

 Now we will see some positive results for the very first classes of Borel sets. We know
(see [L1]) that if $X$ and $A$ are Polish spaces, then every Borel subset of
$X\times A$ with vertical sections in $\boraone(A)$ is 
$\big(\borel (X)\times\boraone (A)\big)_\sigma$. 

\begin{prop} If $A$ has a countable basis, then every $\bormone (X\times A)$ with vertical sections in 
$\boraone (A)$ is $\big(\bormone (X)\times \boraone (A)\big)_\sigma$. If moreover $A$ is 0-dimensional, then every $D_2(\boraone )(X\times A)$ with vertical sections in $\boraone  (A)$ is 
$\big(\bormone (X)\times \borone (A)\big)_\sigma$.\end{prop}

\vfill\eject

\noindent\bf Proof.\rm\ Let $F$ be a closed subset of $X\times A$ with
vertical sections in $\boraone(A)$. As in the proof of Proposition 23, $F = \bigcup_{n}
X_n\times A_n$, where $(A_n)$ is a basis for the topology of $A$. But as $F$ is closed, we also
have $F =
\bigcup_{n} \overline{X_n}\times A_n \in \big(\bormone (X)\times
\boraone(A)\big)_\sigma$.\bigskip

 If $A$ is a 0-dimensional space, let $U$ (respectively $F$) be an open (respectively closed)
subset of $X\times A$ such that $U\cap F$ has vertical sections in $\boraone(A)$; then
$U=\bigcup_n U_n$, where 
$$U_n\in\bormone (X)\times \borone (A).$$ 
For each $x\in X$, we have 
$$(U\cap F)_x = U_x\cap F_x = \bigcup_n ~(U_n)_x\cap F_x=\bigcup_n ~(U_n\cap F)_x.$$
Moreover, $(U_n\cap F)_x = (U_n)_x\cap (U\cap F)_x$ is $\boraone(A)$, so $U_n\cap F$ is 
$\bormone (X\times A)$ with vertical sections in $\boraone(A)$. By the previous facts,
$U_n\cap F\in\big(\bormone (X)\times\borone (A)\big)_\sigma$ and $U\cap F = \bigcup_n U_n\cap
F$ too.$\hfill\square$\bigskip

\begin{prop} There exists a
$\check D_2(\boraone)$ subset of
$2^\omega\times2^\omega$ with sections in $\borone (2^\omega)$ which is not in 
$\big(\bormone (2^\omega)\times\boraone(2^\omega)\big)_\sigma$.\end{prop} 

\noindent\bf Proof.\rm\ This result is a consequence of Proposition 
24. But we can find here a simpler counter-example. We will use it later. Let
$\psi:\omega\rightarrow P_f$ be a bijective map, and 
$${E:=(2^\omega\times\{ 0^\omega\})\cup\bigcup_p ~
(2^\omega\setminus\{\psi(p)\}\times N_{0^p1})}.$$ 
Then $E$ is the union of a closed set and of an open set, so it is 
$\check D_2(\boraone)(2^\omega\times2^\omega)$. If $\alpha\notin P_f$ (respectively $\alpha =\psi (p)$), then we have $E_\alpha = 2^\omega$ (respectively $2^\omega
\setminus N_{0^p1}$); so $E$ has vertical sections in $\borone (2^\omega )$. If ${E =
\bigcup_n F_n\times U_n}$, then we have ${E^{0^\omega} = 2^\omega =
\bigcup_{n/0^\omega\in U_n} F_n}$. By Baire's theorem, there exists $s\in 2^{<\omega}$ and an
integer $n_0$ such that $0^\omega\in U_{n_0}$ and $N_s\subseteq F_{n_0}$. From some integer 
$p_0$ on, we have $N_{0^p1}\subseteq U_{n_0}$. As $P_f$ is dense, there exists
$p\geq p_0$ such that $\psi (p)\in N_s$. We have 
$$\big(\psi (p),0^p10^\omega\big)\in (N_s\times N_{0^p1})\setminus E\subseteq 
(F_{n_0}\times U_{n_0})\setminus E\subseteq E\setminus E.$$
This finishes the proof.$\hfill\square$\bigskip

 Now we will show that the example in Corollary 16 is in some way optimal. Recall 
that the Wadge hierarchy (the inclusion of classes obtained by continuous pre-images of a
Borel subset of $\omega^\omega$; see [LSR]) is finer than that of Baire. The beginning of
this hierarchy is the following:
$$\matrix{\{\emptyset\} & &
\boraone & & D_2(\boraone) &  &
\boratwo\cr & \borone  & & {\boraone}^+ & & \cdots & \cr
\check\{\emptyset\} & &
\bormone  & & \check D_2(\boraone) & & \bormtwo}$$
The class ${\boraone}^+$ is defined as follows: ${\boraone}^+:= \{(U\cap
O)\cup (F\setminus O)~/~U\in\boraone, F\in\bormone ,~O\in\borone\}$.

\vfill\eject

\begin{prop} Let $A$ be a metrizable 
compact space, 
$B\subseteq X\times A$ with vertical (resp., horizontal) sections in $\borone (A)$ $($resp., 
${\boraone}^+(X))$. Then 
${B\in\big(\bormone (X)\times{\cal P}(A)\big)_\sigma}$. In particular, if $Y=2$ and $A$
is made of characteristic functions of ${\boraone}^+(X)$, then conditions 
of Proposition 22
are satisfied and $\phi$ is Baire class one.\end{prop}

\noindent\bf Proof.\rm\ For $f\in A$, we have $B^f=(U^f\cap
O^f)\cup (F^f\setminus O^f)$. We set
$$B_1 := \{(x,f)\in X\times A~/~x\in U^f
\cap O^f\}\mbox{,}~~~~~~ B_2 := \{(x,f)\in X\times A~/~x\in F^f\setminus O^f\}.$$ 
Therefore we have $B = B_1\cup B_2$. Let $(X_n)$ be a basis for the topology of $X$. We have 
$$B_1 =\bigcup_n X_n\times\{f\in A~/~X_n\subseteq O^f\cap U^f\}.$$ 
Thus $B_1\in\big(\bormone (X)\times{\cal P}(A)\big)_\sigma$. In the same way, 
${\{(x,f)\in X\times A~/~x\notin O^f\} =\bigcup_n X_n\times E_n}$, where 
$E_n := \{f\in A~/~X_n\cap O^f=\emptyset\}$. Let us enumerate $\borone (A) :=
\{O_m~/~m\in \xi\}$, where $\xi\in\omega +1$. We have $B_2 = \bigcup_{n,m}~\{x\in
X_n~/~O_m\cap E_n\subseteq B_x\}\times(O_m\cap E_n)$. It is enough to see that $\{x\in
X_n/O_m\cap E_n\subseteq B_x\}\! \in\! \bormone (X_n)$. Let $(f^n_p)_p$ be a dense sequence
of $E_n$. If $x\in X_n$, then 
$$\begin{array}{ll}
O_m\cap E_n\subseteq B_x\!\!\!\! 
& \Leftrightarrow\forall p\in\omega~~f^n_p\notin O_m\cap E_n~~\mbox{or}~~x\in B^{f^n_p}\cr 
& \Leftrightarrow\forall p\in\omega~~f^n_p\notin O_m\cap E_n~~\mbox{or}~~x\in F^{f^n_p}\setminus
O^{f^n_p}.
\end{array}$$ 
Therefore, $B_2\in \big(\bormone (X)\times{\cal P}(A)\big)_\sigma$ and $B$ too.$\hfill\square$\bigskip

\begin{prop} Assume that $X$ and $A$ are Polish spaces, that
$Y=2$, and that $A$ is made of characteristic functions of $D_2(\boraone)(X)$. Then 
$\phi^{-1}(\{ 1\})\in \big(\bormone (X)\times{\cal P}(A)\big)_\sigma$.\end{prop}

\noindent\bf Proof.\rm\ As $\phi$ is Baire class two,
$\phi^{-1}(\{ 1\})$ is $\borthree (X\times A)$ with horizontal sections in $D_2(\boraone)(X)$.
So there exists a finer Polish topology $\tau$ on $A$ and some open subsets $U_0$ and 
$U_1$ of
$X\times [A,\tauÊ]$ such that $\phi^{-1}(\{ 1\}) = U_1\setminus U_0$. The reader should see 
[L1] and [L2] to check this point (it is showed for Borel sets with sections
in $\boraxi$ in [L1]; we do the same thing here, using the fact, showed in [L2], that two disjoint $\Ana$ which can be separated by a $D_2(\boraone)$ set can be separated by
a $D_2(\boraone\cap\Borel )$ set). Let $(A_q)$ (resp., $(X_n)$) be a basis for the topology of
$A$ (resp., $X$). Let 
${E_n := \{f\in A~/~X_n\subseteq U_1^f\}}$. There exists $F^n_l\in\bormone (X)$ such that ${U_1
=\bigcup_n ~X_n\times E_n =\bigcup_{n,l}~ F^n_l\times E_n}$. We set 
$$\begin{array}{ll}
F^{n,l}\!\!\!\! 
& := [F^n_l\times E_n]\cap\phi^{-1}(\{ 1\}) =[F^n_l\times E_n]\setminus U_0\cr 
& = \bigcup_q~\{x\in F^n_l~/~A_q\cap E_n\subseteq \phi^{-1}(\{ 1\})_x\}\times (A_q\cap E_n).
\end{array}$$
This is a closed subset of $F^n_l\times [E_n,\tau\lceil E_n]$, and union of the
$F^{n,l}$'s is $\phi^{-1}(\{ 1\})$. So we have ${\phi^{-1}(\{ 1\}) =
\bigcup_{n,l,q}~\overline{\{x\in F^n_l~/~A_q\cap E_n\subseteq \phi^{-1}(\{
1\})_x\}}\times (A_q\cap E_n)\in\big(\bormone (X)\times{\cal
P}(A)\big)_\sigma}$.$\hfill\square$\bigskip 

 These last two propositions show that the example in Corollary 16 is optimal.
In this example, one has $\phi^{-1}(\{ 0\})\notin \boratwo(X\times A)\cup
\big(\bormone (X)\times{\cal P}(A)\big)_\sigma$.

\vfill\eject

\noindent\bf (D) The~case~of~Banach~spaces.\rm\bigskip

 The reader should see [DS] for basic facts about Banach spaces. Let $E$ be a Banach
space, $X := [B_{E^*},w^*]$, $Y := \mathbb{R}$ and 
${A := \{G\lceil X~/~G\in B_{E^{**}}\}}$. If $E$ is separable, then $X$ is a metrizable compact space. If 
moreover $E$ contains no copy of $l_1$, Odell and Rosenthal's theorem gives, 
for every $G\in E^{**}$, a sequence $(e_p)$ of $E$ such that 
$f(e_p)\rightarrow G(f)$ for each $f\in E^*$ (see [OR]). Let 
$i:E\rightarrow E^{**}$ be the  canonical  map, and $G_p := i(e_p)$. Then $(G_p)$ pointwise 
tends to $G$. By definition of the weak* topology, we have $i(e)\lceil X\in {\cal
C}(X,Y)$ for each $e\in E$, thus $G\lceil X$ is the pointwise limit of a sequence of continuous
functions. Therefore, $G\lceil X\in {\cal B}_1(X,Y)$ (see page 386 in [Ku]). We set 
$$\Phi :\left\{\!\!
\begin{array}{ll} 
[B_{E^{**}},w^*] 
& \!\!\!\!\rightarrow [{\cal B}_1(X,Y),\mbox{p.c.}]\cr
G 
& \!\!\!\!\mapsto G\lceil X
\end{array}
\right.$$ 
By definition of weak* topology, $\Phi$ is continuous, and its range is $A$. So $A$ is a 
compact space because $\Phi$'s domain is a compact space.\bigskip

 If $E^*$ is separable, then $E$ is separable and $E$ contains no copy of $l_1$. Indeed, if 
$\phi$ was an embedding of $l_1$ into $E$, then the adjoint map $\phi^*:E^*\rightarrow l_1^*$ of
$\phi$ would be onto, by the Hahn-Banach theorem. But $l_\infty \simeq l_1^*$ would be
separable, which is absurd. The domain of $\Phi$ is a metrizable compact space, thus it is a
Polish space. Therefore,
$A$ is an  analytic compact space. So it is metrizable (see Corollary 2 page 77 of Chapter 9
in [Bo2]). In particular, every point of $A$ is
$G_\delta$.  Conversely, if $E^*$ is not separable, then $\{0_{E^{**}}\}$ is not a 
$G_\delta$ subset of $B_{E^{**}}$. Indeed, if $(x_p)\subseteq E^*$, closed subspace spanned 
by $\{x_p~/~p\in\omega\}$ is not $E^*$ (see page 5 in [B1]), and we use the
Hahn-Banach theorem. Thus 
$\{0_{E^{**}}\lceil X\}$ is not a $G_\delta$ subset of $A$, because $\Phi$ is continuous. 
So the following are equivalent: $E^*$ is separable,  $A$ is metrizable, and every point of
$A$  is $G_\delta$.\bigskip

 Assume that $E^{**}$ is separable. Then $E^*$ is separable, and $A$ is 
 uniformly recoverable. Indeed,
$A\subseteq {\cal C}([B_{E^*},\Vert .\Vert],Y)$, and the following map is continuous:
$$\Phi':\left\{\!\!
\begin{array}{ll}
[E^{**},\Vert .\Vert] 
& \!\!\!\!\rightarrow [{\cal C}([B_{E^*},\Vert .\Vert],Y),\Vert .\Vert_{\infty}]\cr 
G 
& \!\!\!\!\mapsto G\lceil B_{E^*}
\end{array}
\right.$$ 
Therefore, $[\Phi'[E^{**}],\Vert .\Vert_{\infty}]$ is a separable metrizable space and contains $A$. 
Then we can apply Theorem 20. But we have a better result:

\begin{thm} Let $E$ be a Banach space, $X := [B_{E^*},w^*]$, ${A :=\{G\lceil X/G\in B_{E^{**}}\}}$, and 
$Y:=\mathbb{R}$. The following statements are equivalent:\smallskip

\noindent (a) $E^*$ is separable.\smallskip

\noindent (b) $A$ is metrizable.\smallskip

\noindent (c) Every singleton of $A$ is $G_\delta$.\smallskip

\noindent (d) $A$ is uniformly recoverable.\end{thm}
  
\noindent\bf Proof.\rm\ We have seen that conditions (a), (b) and (c) are equivalent. So let us show that 
(a) $\Rightarrow$ (d). We have seen that
$X$ and $A$ are metrizable compact spaces, and that $A\subseteq {\cal B}_1(X,Y)$. Thus we
can apply Proposition 22, and it is enough to check that condition
(2) is satisfied.

\vfill\eject

 The finer topology is the norm topology. Let us check that it is made
of $\boratwo(X)$. We have $\Vert f-f_0\Vert <\varepsilon ~\Leftrightarrow 
~\exists n~~\forall x\in
B_E~~\vert f(x)-f_0(x)\vert\leq\varepsilon -2^{-n}$.\bigskip

\noindent (d) $\Rightarrow$ (c) Let $G\in A$. Then $\{ G\}=A\cap\bigcap_{p,q}\{g\in
\mathbb{R}^X~/~\vert g(x_p)-G(x_p)\vert <2^{-q}\}$. Thus $\{ G\}$ is $\bormtwo(A)$.$\hfill\square$\bigskip
 
  So we get a characterization of the separability of the dual space of an arbitrary Banach
space. Notice that the equivalence between metrizability of the the compact space and the fact that
each of its point is $G_\delta$ is not true for an arbitrary compact set of Baire
class one functions (because of the ``split interval").\bigskip

 This example of Banach spaces also shows that the converse of Theorem 
 20 is false. Indeed, we
set $X:=[B_{l_1},\sigma(l_1\!,c_0)]$, ${A :=\{G\lceil X/G\in B_{l_\infty}\}}$, and $Y:=\mathbb{R}$.
By Theorem 30, $A$ is uniformly recoverable, since $l_1$ is
separable. But since 
$X$ is compact, compact open topology on $A$ is the uniform convergence topology. If
$A$ was separable for compact open topology,
$l_\infty$ would be separable, which is absurd. Indeed, if $(G_n)\subseteq
B_{l_\infty}$ is such that $\{G_n\lceil X~/~n\in\omega\}$ is a dense subset of $A$ for  
the uniform convergence topology, we can easily check that $\{q.G_n~/~q\in{\mathbb{Q}}_+~{\rm
and}~n\in\omega\}$ is dense in $l_\infty$. Notice that this gives an example of a metrizable 
compact space for the pointwise convergence topology which is not separable for the compact
open topology.\bigskip

 Finally, notice that the map $\phi$ is Baire class one if $E^*$ is separable. Indeed, it is the
composition of the identity map from $X\times A$ into $[X,\Vert .\Vert ]\times A$ (which is
Baire class one), and of the map which associates $G(f)$ to 
$(f,G)\in [X,\Vert .\Vert ]\times A$ (which is continuous).\bigskip

\noindent\bf (E) The~notion~of~an~equi-Baire~class~one~set~of~functions.\rm\bigskip

 We will give a characterization of Baire class one functions which lightly improves, in the 
sense (a) $\Rightarrow$ (b) of Corollary 33 below, the one we can find in [LTZ].
 
\begin{defi} Let $X$ and $Y$ be metric spaces, and
$A\subseteq Y^X$. Then $A$ is $equi\mbox{-}Baire~class~one$ (EBC1) if, 
for each $\varepsilon >0$, there exists $\delta (\varepsilon )\in {\cal B}_1(X,\mathbb{R}^*_+)$ such that 
$$d_X(x,x')< \mbox{min}\big(\delta (\varepsilon )(x),\delta (\varepsilon )(x')\big)\Rightarrow
\forall f\in A~~d_Y\big(f(x),f(x')\big)<\varepsilon .$$\end{defi}

\begin{prop} Let $X$ and $Y$ be metric spaces. Assume that $X$
is separable, that all the closed subsets of $X$ are Baire spaces, and that $A\subseteq Y^X$.
The following conditions are equivalent:\medskip

\noindent (1) $A$ is EBC1.\smallskip

\noindent (2) For each $\varepsilon >0$, there exists a sequence $(G^\varepsilon _m)_m\subseteq\bormone  (X)$, whose
union is $X$, such that for each $f\in A$ and for each integer $m$, $\mbox{diam}(f[G^\varepsilon _m])<\varepsilon $.\smallskip

\noindent (3) There exists a finer metrizable separable topology on $X$, made of
$\boratwo(X)$, making $A$ equicontinuous.\smallskip

\noindent (4) Every nonempty closed subset $F$ of $X$ contains a point $x$ such that $\{f_{\vert
F}~/~f\in A\}$ is equicontinuous at $x$.\end{prop}

\vfill\eject

\noindent\bf Proof.\rm\ (1) $\Rightarrow$ (2) We set, for $n$ integer, 
$H_n:=\{x\in X~/~\delta (\varepsilon )(x)>2^{-n}\}$. As $\delta (\varepsilon )$ is Baire class one, there 
exists ${(F^n_q)_q\subseteq \bormone (X)}$ such that $H_n =
\bigcup_q F^n_q$. We construct, for $\xi <\omega_1$, open subsets $U_\xi$ of $X$, and integers
$n_\xi$ and $q_\xi$ satisfying ${\bigcup_{\eta <\xi} U_\eta\not=X\Rightarrow\emptyset\not=
U_\xi\setminus (\bigcup_{\eta <\xi} U_\eta)\subseteq F^{n_\xi}_{q_\xi}}$. It is
clearly possible since $X = \bigcup_{n,q} F^n_q$ and $X\setminus (\bigcup_{\eta <\xi}
U_\eta)$ is a Baire space. As $X$ has a countable basis, there exists $\gamma<\omega_1$ such
that $\bigcup_{\xi <\omega_1} U_\xi = \bigcup_{\xi \leq\gamma} U_\xi$. In particular we have
${U_{\gamma+1}\subseteq \bigcup_{\xi
\leq\gamma} U_\xi}$, thus $X=\bigcup_{\xi \leq\gamma} U_\xi=\bigcup_{\xi \leq\gamma,\mbox{disj.}}
U_\xi\setminus (\bigcup_{\eta<\xi} U_\eta)$. Let $(x^\xi_q)_q\subseteq X$
satisfying $U_\xi\subseteq \bigcup_q B(x^\xi_q,2^{-n_\xi-1}[$. Let ${G^\varepsilon _{q,\xi} :=
B(x^\xi_q,2^{-n_\xi-1}[\cap U_\xi\setminus (\bigcup_{\eta<\xi} U_\eta )}$. Then 
$G^\varepsilon _{q,\xi}\in \boratwo (X)$ and $X$ is the union of
the sequence $(G^\varepsilon _{q,\xi})_{q,\xi\leq\gamma}$. If ${x,x'\in G^\varepsilon _{q,\xi}}$, then we have 
$d_X(x,x')<2^{-n_\xi}<\mbox{min}\big(\delta (\varepsilon )(x),\delta (\varepsilon )(x')\big)$. Thus 
$${d_Y\big(f(x),f(x')\big)\! <\!\varepsilon }$$ 
for each function $f\in A$. It remains to write the $(G^\varepsilon _{q,\xi})_{q,\xi\leq\gamma}$'s as countable unions of closed sets. So that we get the sequence $(G^\varepsilon _m)_m$.\bigskip

\noindent (2) $\Rightarrow$ (3) Let us take a look at the proof of the implication (1) $\Rightarrow$ 
(2) in Proposition 22. There exists a finer metrizable separable topology on $X$, made of
$\boratwo(X)$, and making $G_m^{2^{-r}}$'s open. This is enough (notice that we do not use
the fact that every closed subset of $X$ is a Baire space to show this implication).\bigskip

\noindent (3) $\Rightarrow$ (4) Let $(X_n)$ be a basis for the finer topology. As $X_n\in\boratwo (X)$, 
$F_n := (F\cap X_n)\setminus\mbox{Int}(F\cap X_n)$ is a meager $\boratwo$ subset of $F$.
Thus $F\setminus (\bigcup_n F_n)$ is a comeager $G_\delta$ subset of $F$. As $F$ is a Baire
space, this $G_\delta$ subset is nonempty. This gives the point $x$ we were looking for.
Indeed, let us fix $\varepsilon >0$. Let $n$ be an integer such that $x\in X_n$ and 
$\mbox{sup}_{f\in A}~\mbox{diam}(f[X_n])<\varepsilon $. Then 
$x\in\mbox{Int}(F\cap X_n)$ and $\mbox{sup}_{f\in A}~\mbox{diam}(f_{\vert F}[\mbox{Int}(F\cap
X_n)])<\varepsilon $.\bigskip

\noindent (4) $\Rightarrow$ (2) Let us fix $\varepsilon >0$. We construct a sequence $(U_\xi)_{\xi <\omega_1}$ of open subsets of $X$ such that $\mbox{sup}_{f\in A}~\mbox{diam}(f_{\vert X\setminus
(\bigcup_{\eta <\xi} U_\eta )} [U_\xi\setminus (\bigcup_{\eta <\xi} U_\eta )])<\varepsilon $ and
$U_\xi\setminus (\bigcup_{\eta <\xi} U_\eta )\not=\emptyset$ if $\bigcup_{\eta <\xi}
U_\eta\not=X$. As in the proof of the implication (1) $\Rightarrow$ (2), there exists
$\gamma <\omega_1$ such that $X=\bigcup_{\xi\leq\gamma} U_\xi$. It remains 
to write the 
$\big(U_\xi\setminus (\bigcup_{\eta <\xi} U_\eta )\big)_{\xi\leq\gamma}$'s as countable unions
of closed sets to get the sequence $(G^\varepsilon _m)_m$.\bigskip

\noindent (2)  $\Rightarrow$ (1) For $x\in X$, we set $m^\varepsilon (x) :=
\mbox{min}\{m\in\omega~/~x\in G^\varepsilon _m\}$, and 
$$\delta (\varepsilon ):\left\{\!\!
\begin{array}{ll}
X 
& \!\!\!\!\rightarrow \mathbb{R}^*_+ \cr 
x 
& \!\!\!\!\!\mapsto d_X\big(x,\bigcup_{r<m^\varepsilon (x)} G^\varepsilon _r\big)
\end{array}
\right.$$
Then $\delta (\varepsilon )$ is Baire classe one since if $A,B>0$, then we have 
$$A<\delta (\varepsilon )(x)<B~\Leftrightarrow ~\exists m~~[x\in G^\varepsilon _m~~\mbox{and}~~\forall r<m~~x\notin G^\varepsilon _r~~\mbox{and}~~A<d_X\big(x,\bigcup_{r<m} G^\varepsilon _r\big)<B].$$
If $d_X(x,x')<\mbox{min}\big(\delta (\varepsilon )(x),\delta (\varepsilon )(x')\big)$, then we have
$x'\notin \bigcup_{r<m^\varepsilon (x)} G^\varepsilon _r$, and conversely. Therefore, $m^\varepsilon (x) = m^\varepsilon (x')$ and
$x,x'\in G^\varepsilon _{m^\varepsilon (x)}$. Thus ${d_Y\big(f(x),f(x')\big)\leq\mbox{diam}(f[G^\varepsilon _{m^\varepsilon (x)}])\! <\!\varepsilon }$, for each function $f\in A$ (notice that we do not 
use the fact that every closed subset of $X$ is a Baire space to show these last two
implications).$\hfill\square$

\vfill\eject

\begin{cor} Let $X$ and $Y$ be metric spaces.
Consider the following statements:\medskip

\noindent (a) $f$ is Baire class one.\smallskip

\noindent (b) $\forall \varepsilon >0~\exists\delta (\varepsilon )\in {\cal B}_1(X,\mathbb{R}^*_+)~~d_X(x,x')<\mbox{min}\big(\delta (\varepsilon )(x),\delta (\varepsilon )(x')\big)~\Rightarrow
~d_Y\big(f(x),f(x')\big)<\varepsilon $.\medskip

\noindent (1) If $Y$ is separable, then (a) implies (b).\smallskip

\noindent (2) If $X$ is separable and if every closed subset of $X$ is a Baire space, then
(b) implies (a).\end{cor}

\noindent\bf Proof.\rm\ To show condition (1), the only thing to notice is the
following. Let $(y_n)\subseteq Y$ satisfying ${Y =
\bigcup_n B(y_n,\varepsilon /2[}$. By condition (a), let ${(F^n_q)_q\subseteq \bormone (X)}$
satisfying ${f^{-1}(B(y_n,\varepsilon /2[) = \bigcup_q F^n_q}$. We enumerate the sequence
$(F^n_q)_{n,q}$, so that we get $(G^\varepsilon _m)_m$. We have ${G^\varepsilon _m\in \bormone  (X)}$,
${X=\bigcup_{m} G^\varepsilon _m}$, and
${\mbox{diam}(f[G^\varepsilon _m])\! <\!\varepsilon }$ for each integer $m$. Then we use the proof of
implication (2) $\Rightarrow$ (1) in Proposition 32.$\hfill\square$\bigskip

\noindent\bf Remark.\rm ~Let $X$ be a Polish space, 
$Y\subseteq\mathbb{R}$, and 
$A\subseteq Y^X$ be a Polish space. We assume that every nonempty
closed subset $F$ of $X$ contains a point of equicontinuity of $\{f_{\vert
F}~/~f\in A\}$. Then, by Proposition 32, $A\subseteq {\cal B}_1(X,Y)$ and 
by Proposition 17, J.
Bourgain's ordinal rank is bounded on $A$. This result is true in a more general context :

\begin{cor} Let $X$ be a metrizable separable space,
$Y\subseteq\mathbb{R}$, 
$A\subseteq Y^X$ and $a<b$ be reals. We assume that every nonempty closed subset $F$ of $X$
contains a point of equicontinuity of $\{f_{\vert F}~/~f\in A\}$. Then $\mbox{sup}\{L(f,a,b)~/~f\in A\}<\omega_1$. In particular, $\mbox{sup}\{L(f)~/~f\in A\}<\omega_1$.\end{cor}

\noindent\bf Proof.\rm\ Using equicontinuity, we construct a sequence
$(U_\xi)_{\xi <\omega_1}$ of open subsets of $X$ satisfying ${\mbox{sup}_{f\in A}~\mbox{diam}
f_{\vert X\setminus (\bigcup_{\eta <\xi} U_\eta )} [U_\xi\setminus (\bigcup_{\eta <\xi}
U_\eta )]<b-a}$ and $U_\xi\setminus (\bigcup_{\eta <\xi} U_\eta)\not=\emptyset$ if 
${\bigcup_{\eta <\xi} U_\eta\not=X}$. As $X$ has a countable basis, there exists 
$\gamma<\omega_1$ such that $X=\bigcup_{\xi\leq\gamma} U_\xi$. Let $G_0 :=
\emptyset$, $G_{\alpha +1} := \cup_{\xi\leq\alpha} U_\xi$ if $\alpha\leq\gamma$,
$G_\lambda := \cup_{\alpha <\lambda} G_\alpha$ if $0<\lambda\leq\gamma$ is a limit ordinal,
and $G_{\gamma +2} := X$. Let us check that, if $f\in A$,
then $(G_\alpha)_{\alpha\leq\gamma +2}\in R(\{f\leq a\},\{f\geq b\})$ (this will be enough). By the
proof of Proposition 32, $f$ is Baire class one. So $\{f\leq a\}$ and $\{f\geq b\}$ are
disjoint $G_\delta$ subsets of $X$. We have $G_\alpha\subseteq \cup_{\xi <\alpha} U_\xi$ if 
$\alpha\leq\gamma +1$, so the sequence $(G_\alpha)_{\alpha\leq\gamma +2}$ is increasing. If 
$\alpha\leq\gamma$ is the successor of 
$\rho$, then ${G_{\alpha +1}\setminus G_\alpha = (\cup_{\xi\leq\alpha} U_\xi )\setminus
(\cup_{\xi\leq\rho} U_\xi) = U_\alpha\setminus (\cup_{\xi <\alpha} U_\xi)}$. So $G_{\alpha
+1}\setminus G_\alpha$ is disjoint from $\{f\leq a\}$ or from $\{f\geq b\}$. If 
$\alpha\leq\gamma$ is a limit ordinal, then 
$${G_{\alpha +1}\setminus G_\alpha =
(\cup_{\xi\leq\alpha} U_\xi )\setminus (\cup_{\xi <\alpha} G_\xi )\subseteq
(\cup_{\xi\leq\alpha} U_\xi )\setminus (\cup_{\xi <\alpha} U_\xi )}$$ 
because $U_\xi\subseteq G_{\xi +1}$ if $\xi <\alpha$. Thus we have the same conclusion. Finally,
$${G_{\gamma +2}\setminus G_{\gamma +1}= X\setminus (\cup_{\xi\leq\gamma} U_\xi )=\emptyset}\mbox{,}$$
and we are done.$\hfill\square$\bigskip

 Now, we will study similar versions of Ascoli's theorems, for Baire class one functions.
A similar version of the first of these three theorems is true:

\begin{prop} If $A$ is EBC1, then $\overline{A}^{\mbox{p.c.}}$ is EBC1.\end{prop}

\vfill\eject

\noindent\bf Proof.\rm\ It is very similar to the classical one. We
set $\delta_{\overline{A}^{\mbox{p.c.}}}(\varepsilon ) := \delta_A(\varepsilon /3)$. Assume that 
$$d_X(x,x')< \mbox{min}\big(\delta_{\overline{A}^{\mbox{p.c.}}} (\varepsilon 
)(x),\delta_{\overline{A}^{\mbox{p.c.}}} (\varepsilon 
)(x')\big)\mbox{,}$$ 
and let $g\in \overline{A}^{\mbox{p.c.}}$. The following set is an open neighborhood of $g$: 
$${O := \{ h\in Y^X~/~d_Y\big(h(x),g(x)\big)<\varepsilon /3~~\mbox{and}~~d_Y\big(h(x'),g(x')\big)<
\varepsilon /3\}}$$ 
(for the pointwise convergence topology). Set $O$ meets $A$ in
$f$. Then we check that 
$$d_Y\big(g(x),g(x')\big)<\varepsilon /3+\varepsilon /3+\varepsilon 
/3=\varepsilon .$$
This finishes the proof.$\hfill\square$\bigskip

 A similar version of the third of Ascoli's theorem is true in one sense:

\begin{prop} Assume that $X$ and $Y$ are separable
metric  spaces, and that $X$ is locally compact. If $A\subseteq {\cal B}_1(X,Y)$, equipped with
the compact open topology, is relatively compact in $Y^X$, then $A$ is EBC1 and $A(x)$ is
relatively compact for each $x\in X$.\end{prop}

\noindent\bf Proof.\rm\ As $X$ is metrizable, $X$ is
paracompact (see Theorem 4, page 51 of Chapter 9 in [Bo2]). By 
Corollary page 71 of Chapter 1 in [Bo1], there exists a locally finite open covering 
$(V_j)_{j\in \omega}$ of $X$ made of relatively compact sets (we use the fact that $X$ is
separable). For $x\in X$, we set 
$$J_x := \{ j\in \omega~/~x\in V_j\}.$$
It is a finite subset of $\omega$. Let
$e (x)\in\omega$ be minimal such that $B(x,2^{-e(x)}[\subseteq \bigcap_{j\in J_x} V_j$.
Notice that $e\in {\cal B}_1(X,\omega)$. Indeed, let $(x^j_q)_q$ be a dense sequence of 
$X\setminus V_j$. We have 
$$e(x)\! =\! p~\Leftrightarrow\left\{\!\!\!\!\!\!
\begin{array}{ll} 
& \exists k~\{\forall j\! >\! k~~x\notin V_j\}~\mbox{and}~x\! \in \! V_k~\mbox{and}~\forall j\!
\leq\!  k~~\{\forall q~~x^j_q\! \notin\!  B(x,2^{-p}[~\mbox{or}~x\notin V_j\}\cr\cr 
& \mbox{and}~~\forall l<p~\exists j\leq k~~\{\exists q~~x^j_q\in B(x,2^{-l}[~\mbox{and}~x\in V_j\}.
\end{array}
\right.$$
$\bullet$ Let us show that $\overline{A}^{\mbox{c.o.}}\subseteq {\cal B}_1(X,Y)$. For $x\in X$, we
let $U_x$ be a relatively compact open neighborhood of $x$. As $X$ is a Lindel\"of space,
$X =
\bigcup_{n} U_{x_n}$; let $K_n := \bigcup_{p\leq n} \overline{U_{x_p}}$. Then $(K_n)$ is an
increasing sequence of compact subsets of $X$ and every compact subset of $X$ is a subset of
one of the $K_n$'s. By Corollary page 20 of Chapter 10 in [Bo2], $Y^X$, equipped with
the compact open topology, is metrizable.\bigskip

 So let $f\in \overline{A}^{\mbox{c.o.}}$. By the previous facts there exists a sequence
$(f_n)\subseteq A$ which tends to $f$, uniformly on each compact subset of $X$. So we have 
$$\forall m\in\omega~~\exists (p^m_n)_n\in\omega^\omega~~\forall x\in
K_m~\forall n\in\omega~~d_Y\big(f(x),f_{p^m_n}(x)\big)<2^{-n}.$$
Therefore, if $F\in\bormone (X)$, then 
$$f^{-1}(F) = \bigcup_{m}~K_m\setminus (\bigcup_{p<m} K_p )\cap\{x\in
K_m~/~\forall n\in\omega~~d_Y\big(F,f_{p^m_n}(x)\big)\leq 2^{-n}\}.$$
We deduce from this that $f^{-1}(F)$ is $G_\delta$, because it is union of countably many
$G_\delta$'s, partitionned by some $\bortwo(X)$. So $f$ is Baire class one
and $\overline{A}^{\mbox{c.o.}}\subseteq {\cal B}_1(X,Y)$.

\vfill\eject

\noindent $\bullet$ Let $f\in \overline{A}^{\mbox{c.o.}}$, $\varepsilon >0$ and $K$ be a compact subset of $X$. We set
$$U(f,\varepsilon ,K) := \{g\in \overline{A}^{\mbox{c.o.}}~/~\forall x\in
K~~d_Y\big(f(x),g(x)\big)<\varepsilon /3\}.$$
Then $U(f,\varepsilon ,K)$ is an open neighborhood of $f$ for the compact open topology, so there exists
an integer $p_{\varepsilon ,K}$ and $(f^{\varepsilon ,K}_i)_{i\leq p_{\varepsilon ,K}}\subseteq \overline{A}^{\mbox{c.o.}}$
such that $\overline{A}^{\mbox{c.o.}} =
\bigcup_{i\leq p_{\varepsilon ,K}} U(f^{\varepsilon ,K}_i,\varepsilon ,K)$, because $\overline{A}^{\mbox{c.o.}}$ is
compact.\bigskip

\noindent $\bullet$ By Corollary 33, if $f\!\in\! \overline{A}^{\mbox{c.o.}}$, then there exists $\delta (f,\varepsilon 
)\! \in\!   {\cal B}_1(X,\mathbb{R}^*_+)$ such that ${d_Y\! \big(f(x),f(x')\big)\! <\!\varepsilon }$ if
$d_X(x,x')< \mbox{min}\big(\delta (f,\varepsilon )(x),\delta (f,\varepsilon )(x')\big)$. We set  
$$\delta (\varepsilon ) :\left\{\!\!
\begin{array}{ll} 
X 
& \!\!\!\!\rightarrow\mathbb{R}^*_+ \cr 
x 
& \!\!\!\!\mapsto\mbox{min}\big(2^{-e(x)},\mbox{min}_{j\in J_x,i\leq p_{\varepsilon /3,{\overline{V_j}}}}
[\delta (f^{\varepsilon /3,{\overline{V_j}}}_i,\varepsilon /3)(x)]\big)
\end{array}
\right.$$
If $d_X(x,x')< \mbox{min}\big(\delta (\varepsilon )(x),\delta
(\varepsilon )(x')\big)$ and $f\in A$, then $d_X(x,x')<2^{-e(x)}$ and ${x'\in \bigcap_{j\in J_x}
V_j}$. Let $j\in J_x$ (so $j\in J_{x'}$) and $i\leq p_{\varepsilon /3,\overline{V_j}}$ be such that $f\in
U(f^{\varepsilon /3,{\overline{V_j}}}_i,\varepsilon /3,\overline{V_j})$. As $x,x'\in\overline{V_j}$ and
$d_X(x,x')<\mbox{min}\big(\delta (f^{\varepsilon /3,{\overline{V_j}}}_i,\varepsilon /3)(x),\delta
(f^{\varepsilon /3,{\overline{V_j}}}_i,\varepsilon /3)(x')\big)$, we have 
${d_Y\big(f(x),f(x')\big)<3.\varepsilon /3=\varepsilon }$. Let us check that $\delta(\varepsilon )$
is Baire class one. If $A,B>0$, $A<\delta (\varepsilon )(x)<B$ is equivalent to 
$$\left\{\!\!\!\!\!\!
\begin{array}{ll} 
& \exists~k~\{\forall~j>k~~x\notin~V_j\}~\mbox{and}~x\in V_k\cr 
& \mbox{and}~\{e(x)>-\mbox{ln}(B)/\mbox{ln}(2)~\mbox{or}~\exists~j\leq k~~x\in V_j~\mbox{and}~\exists~i\leq p_{\varepsilon /3,{\overline{V_j}}}~~\delta (f^{\varepsilon /3,{\overline{V_j}}}_i,\varepsilon /3)(x)<B\}\cr 
& \mbox{and}~\{e(x)<-\mbox{ln}(A)/\mbox{ln}(2)~\mbox{and}~\forall~j\leq k~~x\notin V_j~\mbox{or}~\forall~i\leq p_{\varepsilon /3,{\overline{V_j}}}~~\delta (f^{\varepsilon /3,{\overline{V_j}}}_i,\varepsilon /3)(x)>A\}.
\end{array}
\right.$$
$\bullet$ The last point comes from the continuity of $\phi (x,.)$, for each $x\in X$; this
implies that $\overline{A}^{\mbox{c.o.}}(x)$ is compact and contains $A(x)$.$\hfill\square$\bigskip

\noindent\bf Counter-example.\rm ~A similar version of the second of Ascoli's theorem is
false, in the sense that there are some metric spaces $X$ and $Y$,
$X$ being compact, and a metrizable compact space 
$${A\!\subseteq\! [{\cal B}_1(X,Y),\mbox{p.c.}]}$$
which is EBC1 and such that, on $A$, the compact open topology (i.e., the uniform
convergence topology) and the pointwise convergence topology are different. Indeed, we
set $X:=[B_{l_1},\sigma(l_1\!,c_0)]$, ${A :=\{G\lceil X/G\in B_{l_\infty}\}}$, and
$Y:=\mathbb{R}$. We have seen that $A$ is not separable for the uniform convergence topology. So this
topology is different on $A$ from that of pointwise convergence. Nevertheless,
$A$ is EBC1. Indeed, the norm topology makes $A$ uniformly equicontinuous, and we just have to
apply Proposition 32. Moreover, $A(x)$ is compact
for each $x\in X$ and $A$ is a closed subset of [$\mathbb{R}^X,\mbox{c.o.}$] (we check it in a standard
way). As $A$ is metrizable and not separable in this space, it is not relatively compact.
Therefore, the converse of Proposition 36 is false in general.

\begin{cor} Assume that $X$ and $Y$ are separable metric spaces
and that $X$ is locally compact. If moreover $A\subseteq {\cal
B}_1(X,Y)$, equipped with the compact open topology, is
relatively compact in $Y^X$, then $A$ is uniformly recoverable.\end{cor}

\noindent\bf Proof.\rm\ By Proposition 13 page 66 of Chapter 1
in [Bo1] and Theorem 1 page 55 of Chapter 9 in [Bo2], we can apply
Propositions 32 and 36, and use Proposition 22.$\hfill\square$

\vfill\eject

\noindent\bf Remarks.\rm ~There is another proof of this corollary. Indeed, as in 
the proof of Proposition 36, $Y^X$, equipped with
the compact open topology, is metrizable and $\overline{A}^{\mbox{c.o.}}\subseteq {\cal B}_1(X,Y)$.
Thus 
$\overline{A}^{\mbox{c.o.}}$ is a metrizable compact space for the compact open topology. Thus it
is separable for this topology. Then we apply Theorem 20.\bigskip

 Let $X$ and $Y$ be separable metric spaces. Assume that every closed subset of $X$
is a Baire space, and that $A\subseteq Y^X$. If $A$ is EBC1, then $A\subseteq {\cal B}_1(X,Y)$ and
the conditions of Proposition 22 are satisfied, by Proposition 32. The converse of this is false.
To see this, we use the example of Proposition 27 : $X :=2^\omega$,
$Y := 2$ et $A := \bigcup_p~ \{1\!\! {\rm I}_{2^\omega\setminus\{\psi (p)\}}\}$. By the proof
of (1) $\Rightarrow$ (2) in Proposition 22, there exists a finer metrizable separable
topology $\tau$ on $2^\omega$, made of $\boratwo (2^\omega )$, and making the $\{\psi (p)\}$'s open, for $p\in\omega$. Thus $\tau$ makes the functions of $A$ continuous. But assume that $\tau'$ is a finer metrizable separable topology on $2^\omega$, made of $\boratwo(2^\omega)$, and makes $A$ equicontinuous. We
would have $P_\infty\notin\boraone([2^\omega,\tau'])$. So we could find $\alpha\in P_\infty$
in the closure of $P_f$ for $\tau'$. If $V$ is a neighborhood of $\alpha$ for
$\tau'$, we could choose $\psi (p)\in V\cap P_f$. We would have ${\vert 1\!\! {\rm
I}_{2^\omega\setminus\{\psi (p)\}}\big(\alpha \big)-1\!\! {\rm I}_{2^\omega\setminus\{\psi
(p)\}}\big(\psi (p)\big)\vert =1}$. But this contradicts the equicontinuity of 
$A$. Then we apply Proposition 32. This also shows the utility of the assumption of relative
compactness in Proposition 36 ($A$ is an infinite countable discrete closed set; so it is not 
 compact, in ${\cal B}_1(2^\omega ,2)$ equipped with the compact open topology).\bigskip\smallskip

\noindent\bf {\Large 5 References.}\rm\bigskip

\noindent [Bo1]\ \ N. Bourbaki,~\it Topologie g\'en\'erale, ch. 1-4,~\rm 
Hermann, 1974

\noindent [Bo2]\ \ N. Bourbaki,~\it Topologie g\'en\'erale, ch. 5-10,~\rm 
Hermann, 1974

\noindent [B1]\ \ J. Bourgain,~\it Some remarks on compact sets of first Baire class,~\rm Bull. S.M.
Belg.~30 (1978), 3-10

\noindent [B2]\ \ J. Bourgain,~\it On convergent sequences of continuous functions,~\rm Bull. S.M.
Belg. S\'er. B~32 (1980), 235-249

\noindent [DE]\ \ U. B. Darji and M. J. Evans,~\it Recovering Baire one functions,~\rm Mathematika~
42 (1995), 43-48

\noindent [DS]\ \ N. Dunford and J. T. Schwartz,~\it Linear operators, Part 1,~\rm 
Interscience publishers, Inc., 1964

\noindent [FV]\ \ C. Freiling and R. W. Vallin,~\it Simultaneous recovery of Baire 
one functions,~\rm Real Analysis Exchange~22 (1) (1996/97), 346-349

\noindent [Ke]\ \ A. S. Kechris,~\it Classical Descriptive Set Theory,~\rm 
Springer-Verlag, 1995

\noindent [KL]\ \ A. S. Kechris and A. Louveau,~\it A classification of Baire class one
functions,~\rm Trans. A.M.S.~318 (1) (1990), 209-236

\noindent [Ku]\ \  K. Kuratowski,~\it Topology,~\rm Vol. 1, Academic Press, 1966

\noindent [LTZ]\ \ P. Y. Lee, W. K. Tang and D. Zhao,~\it An equivalent definition of the
functions of the first Baire class,~\rm Proc. A.M.S.~129 (8) (2000), 2273-2275

\noindent [L1]\ \ A. Louveau,~\it A separation theorem for $\Ana$ sets,~\rm Trans. A. M. S.~260 (1980), 363-378

\noindent [L2]\ \ A. Louveau,~\it Some results in the Wadge hierarchy of Borel sets,~\rm Cabal Sem. 79-81 (A. S. Kechris, D. A. Mauldin, Y. N. Moschovakis, eds), Lect. Notes in Math. 1019 Springer-Verlag (1983), 28-55

\noindent [LSR]\ \ A. Louveau and J. Saint Raymond,~\it The strength of Borel Wadge determinacy,~\rm Cabal Sem. 81-85, Lect. Notes in Math.~1333 Springer-Verlag (1988), 1-30

\noindent [M]\ \ Y. N. Moschovakis,~\it Descriptive set theory,~\rm North-Holland, 1980

\noindent [OR]\ \ E. Odell and H. P. Rosenthal,~\it A double-dual characterization of separable Banach
spaces containing $l_1$,~\rm Israel J. Math.~20 (1975), 375-384

\noindent [T]\ \ S. Todor\v cevi\'c,~\it Compact subsets of the first Baire class,~\rm 
J.A.M.S.~12 (4) (1999), 1179-1212

\vfill\eject

\noindent\bf Acknowledgements.\rm ~I would like to thank A. Louveau
who gave me  a copy of the paper [DE] to give a talk in the descriptive set theory seminar
of the University  Paris 6. This was the origin of this work. I also thank G. Debs, who organizes
with A. Louveau  this seminar, for his questions and his interest for this work. Finally, I
thank G. Godefroy for  having asked the question of the uniformity for compact sets of Baire class
one functions, especially in the case of Banach spaces.\bigskip

\end{document}